\documentclass[12pt,a4paper,twoside, reqno]{amsart}
\usepackage{amsfonts, amsthm, amsmath, amssymb}
\usepackage{mathtools}
\usepackage{mathrsfs,amsmath}
\usepackage{hyperref}
\hypersetup{colorlinks=false}

\usepackage[margin=2.7cm]{geometry}

\usepackage{helvet}

\usepackage{amsthm}
\newtheorem{theorem}{Theorem}
\newtheorem{lemma}{Lemma}[section]
\newtheorem{remark}{Remark}

\newtheorem{corollary}{Corollary}

\begin{document}
\author{Mohd Harun,\,\,Saurabh Kumar Singh}	

\address{ Saurabh Kumar Singh \newline {\em Department of Mathematics and Statistics, Indian Institute of Technology, Kanpur, India; \newline  Email: skumar.bhu12@gmail.com
} }
\address{Mohd Harun \newline {\em Department of Science and Mathematics, Indian Institute of Information Technology, Guwahati, India; \newline  Email: harunmalikjmi@gmail.com
} }

\date{\today}

\begin{abstract}  
In this article, we use the circle method to establish non-trivial bounds
for certain summatory functions and a bilinear sum involving the Fourier coefficients
of a general Hecke–Maass cusp form on  $SL(3, \mathbb{Z})$.
\end{abstract}

\title{Bilinear sum involving $GL(3)$ Fourier coefficients}
	
	\maketitle
\section{Introduction}
Let $X \ge 1$ be a real parameter, and let 
$\mathcal{P}(x_1,x_2,\dots,x_s) \in \mathbb{Z}[x_1,x_2,\dots,x_s]$ 
be a polynomial in $s$ variables with integer coefficients. Define the set
\[
R_{X,\mathcal{P}} 
:= \left\{ \mathcal{P}(x_1,x_2,\dots,x_s) \in \mathbb{N} \;:\; 
1 \le x_i \le X,\; x_i \in \mathbb{N},\; 1 \le i \le s \right\}.
\]
For an arithmetic function $\mathcal{A} : \mathbb{N} \to \mathbb{C}$, define the associated summatory function over the values of $\mathcal{P}$ by
\[
\mathcal{B}(\mathcal{A},X,\mathcal{P})
:= \sum_{1 \le n_1,\dots,n_s \le X} 
\mathcal{A}\big(\mathcal{P}(n_1,n_2,\dots,n_s)\big).
\tag{1}
\]
Equivalently, introducing the representation function
\[
r_{\mathcal{P}}(n) 
:= \#\left\{ (n_1,\dots,n_s) \in [1,X]^s \cap \mathbb{N}^s \;:\; 
\mathcal{P}(n_1,\dots,n_s)=n \right\},
\]
we may write
\[
\mathcal{B}(\mathcal{A},X,\mathcal{P})
= \sum_{n \in R_{X,\mathcal{P}}} r_{\mathcal{P}}(n)\,\mathcal{A}(n).
\tag{2}
\]
In the realm of Analytic Number Theory, estimating the non-trivial bound and understanding the asymptotic formula for the different summatory functions $\mathfrak{B}(\mathcal{A}, X, \mathcal{P})$ is one of the fundamental problems. The complexity of this problem is contingent upon the sparseness of the set $R_{X, \mathcal{P}}$, and the nature of the arithmetical function $\mathcal{A}$. Specifically, when $\mathcal{A}(n)$ is taken to be the generalized divisor function $d_{k}(n)$, which are the normalized Fourier coefficients of simplest $GL(k)$ $L$-function $\zeta^k(s)$, the corresponding summatory function $\mathfrak{B}(d_k, X, \mathcal{P})$ has been studied by numerous mathematicians for different values of $k$ and polynomials $\mathcal{P}(x_1, x_2, ..., x_s)$. As a brief literature survey for such type of results, one can go through the following progresses:
\begin{itemize}
    \item For $k=2$ and $\mathcal{P}(x_1) = x^2_1+a$, see Ingham \cite{R11}, Scourfield \cite{R10}, Hooley \cite{R8}.
    \item For $k= 2$ and $\mathcal{P}(x_1)$ is an irreducible polynomial of any degree, see P. Erd\"{o}s \cite{R7}.
    \item For $k =2$ and $\mathcal{P}(x_1, x_2) = x^2_1 + x^2_1$, see Yu \cite{R12}.
    \item For $k =2$ and $\mathcal{P}(x_1, x_2, x_3) = x^2_1 + x^2_2 + x^2_3$, see Calderón and de Velasco \cite{R13}, Guo and Zhai \cite{R14}.
    \item For $k =3$ and $\mathcal{P}(x_1, x_2) = x^2_1 + x^6_2$, see Friedlander and Iwaniec \cite{R15}.
    \item For $k =3$ and $\mathcal{P}(x_1, x_2, x_3) = x^2_1+ x^2_2 + x^2_3$, see Sun and Zhang \cite{R16}.
\end{itemize}
Understanding the asymptotic formula for $\mathfrak{B}(d_k, X, \mathcal{P})$ with $k \geq 3$ and the polynomial $\mathcal{P}(x_1, x_2,..., x_s)$ is of degree $\geq 3$, is a difficult problem and has not been much explored so far. Motivated by the above results, researchers have also explored the analogous results in a more general setup when the particular functions $d_2(n),  d_3(n)$ are replaced by the Fourier coefficient of general $GL(k)$ $L$-functions.

Consider the arithmetical function $\mathcal{A}(n) = \lambda_{f}(n)$, representing the $n$-th Fourier coefficient of a holomorphic cusp form $f \in S(k_1, \Gamma_0(L), \chi)$ with weight $k_1\geq 4$ and character $\chi$. Here, $\Gamma_0(L)$ denotes the congruence subgroup of $\Gamma \,=\,SL(2, \mathbb{Z})$ with level $L \in \mathbb{N}$. For $\mathcal{P}(x_1) = x^2_1 + ax_1 + b$, Blomer \cite{Blo} established the following asymptotic formula:
\begin{align*}
 \mathfrak{B}(\lambda_{f}, X, \mathcal{P}) = \sum_{n \leq X} \lambda_{f}(\mathcal{P}(n)) = c(f, \mathcal{P}) X + O_{f, \mathcal{P}, \epsilon}(X^{6/7+\epsilon}),
\end{align*}
where $c(f, \mathcal{P}) \in \mathbb{C}$ is a constant that vanishes in many cases, although not universally. Subsequent progress in this context was made by Templier and Tsimerman \cite{TT} and Acharya \cite{Ach}.
\vspace{0.2cm}

As one of the main results of this paper, we prove a general result in the same context. Let $\mathcal{A}(n) = A_{\pi}(1, n)$ denote the $n$-th Fourier coefficients of a Hecke-Maass cusp form $\pi$ for $SL(3, \mathbb{Z})$ and $\mathcal{P}(x_1, x_2) = x_1 + x^k_2\,$ for any natural  number $k > 1$. We prove cancellations in the summatory function
\begin{align*}
    \sum_{n \leq X}\,\sum_{m \leq Y}\,A_{\pi}(1, n+m^k),
\end{align*}
where $Y \asymp X^{1/k}$. To achieve our result, we will prove two general theorems that are inherently interesting and have appealing special cases that connect to existing noteworthy results. This is an extension of our earlier work \cite{HC} for $GL(2)$. Utilizing a variant of the delta method introduced by Duke, Friedlander, and Iwaniec, we prove the following theorem.  

\begin{theorem}\label{th1}
Let $M, N \geq 1$ be variables and $\alpha(n), \beta(n)$ be any sequences of complex numbers. Let $A_{\pi}(1,n)$ be the normalized $n$-th Fourier coefficients of a Hecke-Maass cusp form $\pi$ for SL$(3,\mathbb{Z})$. For $M < N^{1-\epsilon}$, we have
\begin{align*}
 S(N,M)  = \, \sum_{1 \leq n \leq N}\,\sum_{1 \leq m \leq M} \, \alpha(n)\,\beta(m)\, A_{\pi}(1,n+m) \ll_{\pi, \epsilon} \,||\alpha||_2 \, ||\beta||_2\,\frac{N^{15/16}}{M}\, (MN)^{1/2 + \epsilon}.   
\end{align*}
This bound for $S(N,M)$ is non-trivial within the range $N^{15/16} < M < N^{1-\epsilon}$. 
\end{theorem}
\vspace{0.3cm}

\noindent In $1993$, Friedlander and Iwaniec studied the GL$(1)$ analog of a similar bilinear sum as in Theorem \ref{th1}. More specifically, let $p$ be a prime and $\chi$ be a non-principal character modulo $p$, Friedlander and Iwaniec \cite{FI} studied the following sum
\begin{align*}
    \sum_{m \in \mathcal{B}_1}\,\sum_{n \in \mathcal{B}_2}\,\chi(m+n),
\end{align*}
where $\mathcal{B}_1$ and $\mathcal{B}_2$ are finite sets of integers. They proved several estimates for this sum, which gives, in particular, a modified proof of the inequalities by Polya-Vinogradov and Burgess. We also proved non-trivial estimates for the similar $GL(2)$ bilinear sums \cite{HC}.
\vspace{0.2cm}

\noindent In $2013$, while studying the GL$(2)$ analog of the classical Titchmarsh divisor problem, Pitt \cite{NP} also derived a non-trivial bound for a similar bilinear sum as in Theorem \ref{th1}. Let $f \in S(k_1, \Gamma)$ and $M, N \geq 1$ be two variables, Pitt considered the following sum,
\begin{align*}
\mathfrak{D}(\alpha, \beta)\,=\,\sum_{n \leq N}\sum_{m \leq M} \alpha(n)\,\beta(m)\,\lambda_{f}(mn-1).   
\end{align*}
He established the bound,
\begin{align*}
 \mathfrak{D}(\alpha, \beta)\, \ll_{f, \epsilon} ||\alpha||_2\, ||\beta||_2 \left\{M^{3/8}N^{1/2} + M^{3/4} N^{3/8}\right\}\,(MN)^{(k_1-1)/2 + \epsilon}. 
\end{align*}
By Theorem \ref{th1}, we are kind of extending our results \cite{HC} and these results of Friedlander, Iwaniec, and Pitt by studying the GL$(3)$ analog of a similar bilinear sum. Now, we will explore some intriguing corollaries derived from Theorem \ref{th1}.
\vspace{0.3cm}

\noindent As a particular case of Theorem \ref{th1}, we get the following results in the context of shifted convolution sums, although we admit that the bound may not be optimal.

\begin{corollary}\label{f8}
Let $M, N$ be defined as in Theorem \ref{th1}. Consider the following cases:
\begin{itemize}
    \item $\alpha(n) = d_k(n)$ and $\beta(m) = 1$.
    \item $\alpha(n) = A_{\pi_1}(1, n)$ and $\beta(m) = 1$.
    \item $\alpha(n) = A_{\pi_1}(1, n)$ and $\beta(m) = A_{\pi_2}(1, m)$
\end{itemize}
where $A_{\pi_1}(1, n)$ and $A_{\pi_2}(1, n)$ represent the normalized $n$-th Fourier coefficients of two not necessarily same $SL(3, \mathbb{Z})$ Hecke-Maass cusp forms $\pi_1$ and $\pi_2$, respectively. In all these cases, using Theorem \ref{th1}, we have
\begin{align*}
\,\frac{1}{M}\,\mathop{\sum_{m\leq M}\,\beta(m)\,\sum_{n\leq N}}\, \alpha(n)\,A_{\pi}(1, n + m) \,\ll_{\pi, \epsilon}\, \frac{N^{15/16}}{M} \,N^{1+\epsilon}.   
\end{align*} 
This bound is non-trivial in the range $N^{15/16} < M < N^{1-\epsilon}$.
\end{corollary}

\noindent This corollary extends the findings of Pitt \cite{NJPT}. Also, it establishes new results related to the non-trivial estimates of general $GL(3) \times GL(3)$ shifted convolution sums in both average and weighted average setups.  
\vspace{0.3cm}

\noindent Let $w \in \mathbb{N}$ and $r_w(n) \,=\, \# \left\{(n_1,n_2,...,n_w) \in \mathbb{Z}^w:\,n^2_1 + n^2_2 +...+ n^2_w = n \right\}$ represents the $n$-th Fourier coefficient of the modular form $\theta^w(z)$. Here, $\theta(z)$ represents the classical Jacobi theta series. It is well-known that the estimates $r_w(n) \asymp n^{w/2 -1}$ hold for $w > 4$. For $1 \leq h\leq X$, Sun \cite{R17} established the following result:
\begin{align*}
 \sum_{1\leq n \leq {X}}\,A_{\pi}(1,n+h)\,r_3(n)\,=\mathop{\sum\, \sum\, \sum}_{1\leq n_1,n_2,n_3 \leq \sqrt{X}}\,A_{\pi}(1,n^2_1 + n^2_2 + n^2_3+h)\ll_{\pi, \epsilon}\, X^{3/2-1/8+\epsilon}.
\end{align*}
The above sum can also be interpreted as the shifted convolution sum of arithmetical functions $A_{\pi}(1,n)$ and $r_3(n)$. As a corollary of Theorem \ref{th1} and in the same context, we can also obtain the cancellations in a shifted convolution sum of $r_w(n)$ and $A_{\pi}(1, n)$ with an average over the shift.\\

\begin{remark}
For any two complex sequences $\alpha(n)$ and $\beta(n)$ satisfying the condition 
\begin{align*}
||\alpha||^2_2 =  \sum_{n \leq X}\,|\alpha(n)|^2 \,\ll_{\epsilon}\, X^{1+\epsilon}, \hspace{0.6cm}\text{and}  \,\,\,\,\,\,||\beta||^2_2 \,\ll_{\epsilon}\,X^{1+\epsilon}, 
\end{align*}
we can obtain the same strength bound as in Corollary \ref{f8}.
\end{remark}
\vspace{0.3cm}

\noindent Now, let's set the background for the next main result of this paper. Let $k \in  \mathbb{N}$ and consider a set $S \subset [X, 2X]$ such that $|S| \asymp X^{1/k}$. For example, one can take $S = \{m \in \mathbb{N}: X \leq m^k \leq 2X\}$. By utilizing the result of Theorem \ref{th1}, we establish cancellations in the following sum:
\begin{align*}
T(X, S) =  \sum_{n \leq X}\,\sum_{m \in S} \, \lambda(n)\, A_{\pi}(1,n+m),
\end{align*}
where $\lambda(n)$ represents any sequence of complex numbers. Trivially, one can get the following bound $$T(X, S) \,\ll_{\epsilon}\,||\lambda||_2\, \sqrt{X}\,|S| X^{\epsilon}.$$ 
We prove the following theorem as the second main result of this paper.

\begin{theorem}\label{th2}
Let $\lambda(n)$ be any complex sequence and $A_{\pi}(1,n)$ be defined as in Theorem \ref{th1}.  Let $X \geq 1$ be a variable, and $S \subset[X,2X]$ such that $|S| \asymp X^{1/k}$ with natural number $k >1 $. We have
\begin{align*}
T(X, S) =  \sum_{X \leq n \leq 2X}\,\sum_{m \in S} \, \lambda(n)\, A_{\pi}(1,n+m) \,\ll_{\pi, \epsilon}\,  ||\lambda||_2\,\, \sqrt{X}\,\,\left(||\beta_{\ell}||_2\right)^{1/2^{\ell+1}}\,  X^{\frac{15}{2^{\ell +5}}- \frac{1}{2k} + \epsilon},
\end{align*}
where $\ell$ is least such that $2^{\ell+1} > k (1-\wp)$ for $\wp > 1/16$ and 
\begin{align*}
 |\beta_{\ell}(m)| =  \mathop{\sum_{m_1 \in S}\,\sum_{m_2 \in S}\,...\,\sum_{m_{2^{\ell+1}} \in S}}_{m \,=\,(m_1-m_2) + (m_3 - m_4) \,+\, ...\, + \,(m_{2^{\ell+1}-1}- m_{2^{\ell+1}}) } \,1. 
 \end{align*}
\end{theorem}
\vspace{0.3cm}

\noindent Let $\mathcal{A}(n) = A_{\pi}(1,n)$ and $\mathcal{P}(x_1, x_2) = ax^2_1 + b x_1 x_2 + cx^2_2$ is a positive definite quadratic form in two variables. Recently, Chanana and Singh \cite{R18} established a non-trivial bound for the corresponding summatory function. Specifically, they proved
\begin{align*}
\mathfrak{B}(A_{\pi}, X, \mathcal{P}) \, = \, \mathop{\sum \sum}_{1 \leq n,m \leq X}\,A_{\pi}(1, an^2 + bmn + cm^2)\,  \ll_{\pi, \epsilon} \, X^{2-1/68 + \epsilon}.
\end{align*} The task becomes comparatively more straightforward when considering the polynomial $\mathcal{P}(x_1,x_2)$ to be homogeneous. However, as a consequence of Theorem \ref{th2}, the paper achieves better results when the terms of polynomials are not homogeneous. More precisely, we get the following result:

\begin{corollary}\label{M19}
Let $\lambda(n) = 1$ for all $n \in \mathbb{N}$ and $S = \{m \in \mathbb{N} : X \leq m^2 \leq 2X \}$. We have,
\begin{align*}
\sum_{X \leq n \leq 2X}\,\sum_{m \in S} \, \lambda(n)\, A_{\pi}(1,n+m) =     
\sum_{X \leq n \leq 2X}\,\sum_{\sqrt{X} \leq m \leq \sqrt{2X}} \, A_{\pi}(1,n+m^2) \,\ll_{\pi, \epsilon}\,  X^{\frac{3}{2}- \frac{1}{32} + \epsilon}.    
\end{align*}
\end{corollary}
\begin{proof}
Since  $S = \{m \in \mathbb{N} : X \leq m^2 \leq 2X \}$, then $|S| \asymp X^{1/2}$. For $k = 2$, we get $\ell = 0$. Using Lemma \ref{beta}, we get $||\beta_0|| \,\ll\,X^{1/2}$. We obtain the desired result using these particular values and estimate in Theorem \ref{th2}.     
\end{proof}
\vspace{0.3cm}

\noindent For the higher degree polynomials in two variables of the form $\mathcal{P}(x_1,x_2) = x_1 + x^k_2$ with $k \geq 3$, we can also get cancellations by using Theorem \ref{th2}. We have the following corollary:
\begin{corollary}\label{M119}
Let $\lambda(n) = 1$ for all $n \in \mathbb{N}$ and $S = \{m \in \mathbb{N} : X \leq m^k \leq 2X \}$ with $k \geq 3$. We have,
\begin{align*}
\sum_{X \leq n \leq 2X}\,\,\,\sum_{{X^{1/k}} \leq m \leq {(2X)^{1/k}}} \, A_{\pi}(1,n+m^k) \,\ll_{\pi, \epsilon}\,  X^{1 +\frac{7}{2^{\ell+5}} +\frac{1}{2k} + \epsilon }.    
\end{align*}
This bound is non-trivial when $\ell \in \mathbb{N}$ is least such that $2^{\ell+1} > \text{max}\left\{k(1-\wp),\, k \log k,\, \frac{7k}{8}\right\}$ for $\wp \geq 1/16$.
\end{corollary}
\begin{proof}
 Using Theorem \ref{th2} and Lemma \ref{beta}, the result is immediate.   
\end{proof}

\begin{remark}
For $\lambda(n) = \Lambda(n),\, \mu(n),\, d_{m}(n)$ in Theorem \ref{th2}, we can obtain the same strength bound as in Corollary \ref{M119}.     
\end{remark}

\begin{remark}
The result of Theorem \ref{th2} is general in the sense that, when we choose the set $S \subset [X, 2X]$ such that the corresponding sequence $\beta_{\ell}(m)$ (defined in Lemma \ref{beta}) satisfies a similar bound condition as in Lemma \ref{beta}, cancellations can be obtain using Theorem \ref{th2} in the corresponding summatory functions involving Fourier coefficients of $GL(3)$-automorphic form.
\end{remark} 
\vspace{0.2cm}

\textbf{Notations:} For  $z \in \mathbb{C}$, we set $e(z):= e^{2 \pi i z}$.  By $A \sim B$ we mean that $B \leq A \leq 2B$.  We use $X \ll Y$ to mean that there is a constant $K >0$ such that $|X| \leq K Y $ and by  $A \asymp B$ we mean that $ A \ll B$ and $B \ll A$.  By the term `savings'  we mean  `initial bound'  / ` final bound'.
\vspace{0.4cm}

\section{Sketch of the proof (Theorem \ref{th1})}\label{sketch}
We aim to get cancellations in the following bilinear sum 
\begin{align}\label{j1}
	S(N,M) = \sum_{1 \leq n \leq N}\,\sum_{1 \leq m \leq M} \alpha(n)\,\beta(m)\, A_{\pi}(1,n+m),
\end{align}
 over the trivial bound
 \begin{align*}
	S(N,M) \ll_{\epsilon} ||\alpha||_2 \, ||\beta||_2 \, (MN)^{1/2+\epsilon}.
\end{align*}
where $M < N^{1-\epsilon}$ and $\epsilon >0$ is a small real number. To go through the main steps of the proof, we consider the generic case, i.e., $n \asymp N$ and $m \asymp M$.
 Upon applying the Cauchy-Schwarz inequality and the suitable change of variables, we can write
\begin{align}\label{L18}
    S(N,M) \ll\, ||\alpha||_2\, ||\beta||_2 \,M^{1/4}\, (\eta(n,m))^{1/4},
\end{align}
where $\eta(N, M)$ is given by
\begin{align*}
 \sum_{n_1 \asymp N}\,\sum_{n_2 \asymp N}\,\, A_{\pi}(1,n_1)\,\overline{A_{\pi}(1,n_2)}\, \sum_{n_3 \asymp N}\, W\left(\frac{n_3-n_1}{M}\right) \,A_{\pi}(1, n_3)\,\overline{A_{\pi}(1, n_2 +n_3-n_1)},
\end{align*}
where $W$ is some bump function that records information about a particular variable change and will help in getting additional savings later. Next, we detect the equation $n_2+n_3-n_1 = r$ with the help of delta symbol $\delta(n_2+n_3-n_1, r)$ and apply the DFI delta method (\ref{a1}), the above expression of $\eta(n,m)$ becomes
\begin{align}\label{L3}
   &\, \,\frac{1}{Q}\,\sum_{q \sim Q}\frac{1}{q}\, \sideset{}{^\star} \sum_{a\bmod q}\,
     \sum_{n_1 \asymp N}\,A_{\pi}(1,n_1) \,e\left(\frac{an_1}{q}\right)\,
    \sum_{n_2 \asymp N}\,\overline{A_{\pi}(1,n_2)}\, e\left(\frac{-an_2}{q}\right)\, \notag\\
    &\times\, \sum_{n_3 \asymp N}^{\infty}\,{A_{\pi}(1,n_3)}\, e\left(\frac{-an_3}{q}\right)\,\,W\left(\frac{n_3-n_1}{M}\right) 
  \, \sum_{r \asymp N}^{\infty}\,\overline{A_{\pi}(1,r)}\, e\left(\frac{ar}{q}\right)\, + O(N^{-2025}).
\end{align}
where $q \sim Q \asymp \sqrt{N}$.
Next, we employ the GL$(3)$ Voronoi summation formula (Lemma \ref{gl3voronoi}) to the sum over $n_1, n_2, n_3$, and $r$. Our derivation reveals that the dual length (conductor/initial length) corresponding to the sum over $n_1$, $n_3$ is supported by $(QN)^3/M^3N$, and the sum over $n_2$, $r$, the dual length is supported on $Q^3/N$. Therefore, roughly we arrive at the following expression for $\eta(N, M)$
 \begin{align*}
     \eta(N, M)= \,&\frac{N^{8/3}}{Q^6}\sum_{q\sim Q} \,\, \sum_{n''_{1}\sim (QN)^3/M^3N}\,  \frac{A_{\pi}(n''_{1}, 1)}{ (n''_{1})^{1/3}} \,\,\sum_{n''_{2} \sim Q^3/N}\,  \frac{A_{\pi}(n''_{2}, 1)}{ (n''_{2})^{1/3}} \\
 &\times\ \sum_{n''_{3} \sim (QN)^3/M^3N}\,  \frac{A_{\pi}(n''_{3}, 1)}{ (n''_{3})^{1/3}} \, \, \sum_{r_{2} \sim Q^3/N}  \frac{A_{\pi}(r_{2}, 1)}{ r^{1/3}_{2}}\,\, \mathcal{C}(...)\, \,\mathcal{I}(...) + O(N^{-2024}).
 \end{align*}
Here, $\mathcal{I}(...)$ represents the integral transform as defined in equation \eqref{L11}, and
\begin{align*}
    \mathcal{C}(...) = \sideset{}{^\star} \sum_{a \bmod q}\, S\left(\overline{a}, \pm n''_{1}; q\right) S\left(-\overline{a}, \pm n''_{2}; q\right) S\left(-\overline{a}, \pm n''_{3}; q\right) S\left(\overline{a}, \pm r_{2}; q\right)  
\end{align*}
 denotes the character sum. It's noteworthy that, upon applying the Voronoi summation formula, we achieve savings ($ = \sqrt{\text{initial length}/\text{conductor}}$) of size ${N}/{(QN/M)^{3/2}}$ in the $n_1$ and $n_3$ lengths, and $N/Q^{3/2}$ in the $n_2$ and $r$ lengths.
\vspace{0.3cm}

\noindent  One main challenge was to achieve the desired square root cancellations in the character sum $\mathcal{C}(...)$. This was the main hurdle for passing from $GL(2)$ to $GL(3)$. Using some nice trick in Section \ref{M26}, involving reducing $\mathcal{C}(...)$ to the prime modulus to convert it into some nice form. Then, by leveraging a result from Fouvry, Kowalski, and Michel (see Lemma \ref{t1}), we obtain the desired bound for $\mathcal{C}(...)$. Furthermore, in Section \ref{M27}, we analyze the $4$-fold integral transform $\mathcal{I}(...)$ and attain the desired bound. 
\vspace{0.3cm}

\noindent Against the trivial bound $N^4$ in \eqref{L3} for $\eta(N, M)$, our cumulative savings thus far amount to:
 \begin{align*}
		\frac{N}{Q^{3/2}}\,\times \frac{N}{Q^{3/2}}\,  \times \frac{N}{(QN/M)^{3/2}}\,  \times \, \frac{N}{(QN/M)^{3/2}}\,\times \sqrt{Q} =   \frac{M^3}{N^{7/4}}. 
	\end{align*}
Now, the bound for $\eta(n,m)$ will become
	\begin{align*}
		\eta(N, M) \ll \frac{N^4}{M^3/N^{7/4}}   =  \frac{N^{23/4}}{M^3}.
	\end{align*}
Hence, by \eqref{L18}, we get
	\begin{align*}
		S(N,M) \ll \, ||\alpha||_2\, ||\beta||_2\,M^{1/4}\,\left(\frac{N^{23/4}}{M^3}\right)^{1/4}.    
	\end{align*}
This bound for $S(N,M)$ is non-trivial iff
	\begin{align*}
		M^{1/4}\,\frac{N^{23/16}}{M^{3/4}} < (MN)^{1/2}\,\,\,\,\,  \Longleftrightarrow\,\,\,\,\,  N^{15/16} < M < N^{1-\epsilon}.
	\end{align*}
\vspace{0.4cm}
	
\section{DFI Delta Method}\label{delta}
In this paper, we will use a version of the delta method due to Duke, Friedlander, and Iwaniec. More specifically, we will use the expansion $(20.157)$ given in Chapter 20 of \cite{R1}. Let $\delta : \mathbb{Z}\to \{0,1\}$ be defined by
\[
\delta(n,m)=
\begin{cases}
	1 &\text{if}\,\,n=m \\
	0 &\text{otherwise}.
\end{cases}
\]
Then for $n,m\in\mathbb{Z}\cap [-2L,2L]$, we have
\begin{equation}\label{a1}
	\delta(n,m)=\frac{1}{Q} \sum_{1 \leq q \leq Q} \frac{1}{q} \, \, \sideset{}{^\star} \sum_{a\bmod q}e\left(\frac{(n-m) a}{q}\right)\int_{\mathbb{R}} \psi(q,x)e\left(\frac{(n-m) x}{q	Q}\right) d x, 
\end{equation} where $Q=2L^{1/2}$. Following properties of the function $\psi(q,u)$ are of our interest (see $(20.158)$ and $(20.159)$ of \cite{R1}, and \cite[Lemma 15]{R2})

\begin{align}\label{Delta}
	&\psi(q,x)=1+h(q,x),\,\,\,\,\text{with}\,\,\,\,h(q,x)=O\left(\frac{Q}{q }\left(\frac{q}{Q}+|x|\right)^A\right),\\
	&\psi(q,x)\ll |x|^{-A},\\
	&x^j \frac{\partial^j}{ \partial x^j} \psi(q, x) \ll \  \min \left\lbrace \frac{Q}{q}, \frac{1}{|x|} \right\rbrace \  \log Q.
\end{align}
\vspace{0.2cm}

\noindent for any $A>1$,  $j \geqslant 1$.  In particular, the second property implies that the effective range of integral in \eqref{a1} is $[-L^{\epsilon}, L^{\epsilon}]$. It also follows that if $q \ll Q^{1- \epsilon}$ and $  x \ll Q^{- \epsilon} $, then  $ \psi(q, x)$ can
be replaced by $1$ at the cost of a negligible error term. If $ q \gg Q^{1- \epsilon}$, then we get $ x^j \frac{\partial^j}{ \partial x^j} \psi(q, x) \ll Q^{\epsilon}$, for any $ j \geqslant 1$. If $ q \ll  Q^{1- \epsilon}$  and $  Q^{- \epsilon} \ll |x| \ll   Q^{ \epsilon}$, then $ x^j \frac{\partial^j}{ \partial x^j} \psi(q, x) \ll Q^{\epsilon}$, for any $ j \geqslant 1$. Hence, we can view $  \psi(q, x)$ as a nice weight function in all cases. Finally, by using Parseval and Cauchy-Schwarz inequality, we can get
\begin{align}\label{l2b}
    \int_{|x| \leq Q^{\epsilon}}\, \left(|\psi(q, x)| + |\psi(q, x)|^2  \right)\,dx\,\ll\, Q^{\epsilon},
\end{align}
which shows that the average size of $\psi(q, x)$ is $1$ in $L^1$ and $L^2$ sense. 
\vspace{0.4cm}

\section{Preliminary}

	\subsection{Background on the Maass forms for $SL(3,\mathbb{Z})$} \label{gl3 maass form} 

Let $\pi$ be a Maass form of type $(\nu_{1}, \nu_{2})$ for $SL(3, \mathbb{Z})$. 
We have the Fourier Whittaker expansion for the Maass form (By Jacquet, Piatetskii-Shapiro, and Shalika) $\pi(z)$ as follows
\begin{equation} \label{four-whi-exp}
	\pi(z) = \sum_{\gamma \in U_{2}\left(\mathbb{Z}\right) \backslash  SL(2,\mathbb{Z})} \sum_{m_{1}=1}^{\infty} \sum_{m_{2} \neq 0} \frac{A_{\pi}(m_{1},m_{2})}{m_{1} |m_{2}|} \, W_{J}\left(M \left( {\begin{array}{cc}
			\gamma &  \\
			& 1 \\
	\end{array} } \right)
	z, \nu, \psi_{1,1}\right),
\end{equation} for all $z \in SL(3,\mathbb{Z})\setminus {\mathfrak{h}}^3$.
Where ${\mathfrak{h}}^3$ is the upper half plane in $\mathbb{R}^3$, $U_{2}(\mathbb{Z})$ denotes the group of upper triangular matrices with integer entries and ones on the diagonal,
$W_{J}\left(z,\nu,\psi_{1,1}\right)$ is the Jacquet-Whittaker function, and $M=\textbf{diag} \left(m_{1}|m_{2}|,m_{1},1\right)$ (cf. Goldfeld \cite{R3}). We introduce the Langlands parameters $({\bf \alpha}_{1},{\bf \alpha}_{2},{\bf \alpha}_{3})$ of $\pi$ are given by
$${\bf \alpha}_{1} = - \nu_{1} - 2 \nu_{2}+1, \, {\bf \alpha}_{2} = - \nu_{1}+ \nu_{2},  \, {\bf \alpha}_{3} = 2 \nu_{1}+ \nu_{2}-1.$$ 
We consider $\pi$ to be the eigenfunction of all the Hecke operators with Fourier coefficients $A_{\pi}(m_1, m_2)$, normalized so that $A_{\pi}(1,1) = 1$.
\vspace{0.2cm}

\noindent For $\psi \in C_{c}^{\infty}(0, \infty)$, we denote $\tilde \psi(s) = \int_{0}^{\infty} \psi(x) x^{s-1} \mathrm{d}x$ to be its Mellin transform.For $\sigma > -1 + \max \{-\Re({\bf \alpha}_{1}), -\Re({\bf \alpha}_{2}), \Re({\bf \alpha}_{3})\}$ and $\ell =0$ and $1$, we define
\begin{equation} \label{A15}
	G_{\ell}(y) :=  \frac{1}{2\pi i} \int_{(\sigma)} (\pi^3 y)^{-s} \, \prod_{i=1}^{3} \frac{\Gamma\left(\frac{1+s+{\bf \alpha}_{i}+ \ell}{2}\right)}{\Gamma\left(\frac{-s-{\bf \alpha}_{i}+ \ell}{2}\right)} \tilde \psi(-s)ds, 
\end{equation}

and set \[G_{\pm}(y) = \frac{1}{2\pi^{3/2}} \big(G_0(y) \mp iG_1(y)\big).\]

With the aid of the above terminology, we state the $GL(3)$-Voronoi summation formula in the following lemma.
\begin{lemma} \label{gl3voronoi}
	Let $\psi (x)$ be a compactly supported smooth function on $(0,\infty)$. Let $A_{\pi}(n,m)$ be the $(n,m)$-th Fourier coefficient of a Maass form $\pi(z)$ for $SL(3,\mathbb{Z})$. Then we have
	\begin{align} \label{GL3-Voro}
		& \sum_{n=1}^{\infty} A_{\pi}(n,m)\, e\left(\frac{an}{q}\right) \psi(n) \\
		\nonumber & =q  \sum_{\pm} \sum_{n_{1}|qm} \sum_{n_{2}=1}^{\infty}  \frac{A_{\pi}(n_{1},n_{2})}{n_{1} n_{2}} S\left(m \overline{a}, \pm n_{2}; mq/n_{1}\right) \, G_{\pm} \left(\frac{n_{1}^2 n_{2}}{q^3 m}\right),
	\end{align} 
	where $(a,q)=1, \overline{a}$ is the multiplicative inverse modulo $q$ and $$S(a,b;q) = \sideset{}{^\star}{\sum}_{x \,\rm mod \, q} e\left(\frac{ax+b\overline{x}}{q}\right), $$
	is the Kloostermann sum.
\end{lemma}
\begin{proof}
	See \cite{SD}.
\end{proof}

The following lemma gives an asymptotic expansion for $G_{\pm}(x)$.
\begin{lemma} \label{GL3oscilation}
	Let $G_{\pm}(x)$ be as above,  and  $\psi(x) \in C_c^{\infty}(X,2X)$. Then for any fixed integer $K \geq 1$ and $xX \gg 1$, we have
	\begin{equation*}
		G_{\pm}(x)=  x \int_{0}^{\infty} \psi(y) \sum_{j=1}^{K} \frac{c_{j}({\pm}) e\left(3 (xy)^{1/3} \right) + d_{j}({\pm}) e\left(-3 (xy)^{1/3} \right)}{\left( xy\right)^{j/3}} \, \mathrm{d} y + O \left((xX)^{\frac{-K+5}{3}}\right),
	\end{equation*}
	where $c_{j}(\pm)$ and $d_{j}(\pm)$ are some  absolute constants depending on $\alpha_{i}$,  $i=1,\, 2,\, 3$.  
\end{lemma}
\begin{proof}
We will get our desired result by using \cite[Lemma $6.1$]{XLI}.
\end{proof}

\begin{remark}\label{rm24}
We note that the main term will correspond to $j = 1$ because the higher values of $j$ will reduce the magnitude of the weight function due to the condition $xX \gg 1$ and size of the terms due to $x^{j/3}$ in the denominator. Therefore, our analysis will focus solely on the $j = 1$ term. 
\end{remark}
\vspace{0.3cm}

\noindent The following lemma gives the Ramanujan bound on average for the $GL(3)$ Fourier coefficients $A_{\pi}(m,n)$.

\begin{lemma} \label{ramanubound}
Let $\pi$ be a Maass cusp for SL$(3, \mathbb{Z})$ of type $\nu = (\nu_1, \nu_2) \in \mathbb{C}^2$ with the normalized $n$-th Fourier coefficients $A_{\pi}(m, n)$ and Langland's parameters $\alpha_j,\,1 \leq j \leq 3$. Then we have
\begin{align*}
\mathop{\sum \sum}_{m^{2} n \leq X} \vert A_{\pi}(m,n)\vert ^{2}\, \ll_{\epsilon} \,C^{\epsilon}_{\pi} X^{1+\epsilon},    
\end{align*}
where $C_{\pi} = \prod_{j =1}^{3}\,(1 + |\alpha_j|)$ is the analytic conductor of $\pi$.
\end{lemma}
\begin{proof}
  See  \cite{XI}.
\end{proof}

\noindent In the next lemma, we will prove the square root cancellation in a sum running over all the square full numbers up to a variable $X$.

\begin{lemma}\label{ll2}
Let $X \geq 1$ be a variable, we have
\begin{align*}
    \mathop{\sum_{n \leq X}}_{n \,\text{is square full}}\,1\,\ll_{\epsilon}\, \sqrt{X} X^{\epsilon}.
\end{align*}
\end{lemma}
\begin{proof}
 Since $n$ is square full, we can write $n = m^2_1 m_2$ such that $m_2 | m_1$, we have
 \begin{align*}
      \mathop{\sum_{n \leq X}}_{n \,\text{is square full}}\,1\, \leq \,  \sum_{m_2 |n}\,\,\, \mathop{\sum_{m_1 \leq \sqrt{\frac{X}{m_2}}}}\,1.
 \end{align*}
 Again write $m_1 = m_2 m_3$, we have
 \begin{align*}
\sum_{m_2 |m_1}\,\,\, \mathop{\sum_{m_1 \leq \sqrt{\frac{X}{m_2}}}}\,1 \,\leq \, \sum_{m_2 |n}\,\,\, \mathop{\sum_{m_3 \leq \frac{1}{m_2}\sqrt{\frac{X}{m_2}}}}\,1 \,\leq \,\sqrt{X}\,\sum_{m_2 |n}\, \frac{1}{m_2}{\frac{1}{\sqrt{m_2}}}  \,\ll_{\epsilon}\, \sqrt{X} X^{\epsilon}. 
 \end{align*}
 This proves our lemma.
\end{proof}
\vspace{0.3cm}

\noindent We will use some results related to the estimation of the exponential integral of the form: 
	\begin{equation} \label{eintegral}
		\mathfrak{I}= \int_a^b g(t) e(f(t)) dt,
	\end{equation} where $f$ and $g$ are  real valued smooth functions on the interval $[a, b]$. 
	
\begin{lemma} \label{db}
		Let $f$ and $g$ be real valued twice differentiable functions and let $f^{\prime \prime} \geq r>0$ or  $f^{\prime \prime} \leq -r <0$, throughout the interval $[a, b]$. Let $g(x)/f^\prime(x)$ is monotonic and $|g(x)| \leq A$. Then we have
\begin{align*}
			|\mathfrak{I}|\, \leq\, \frac{8A}{\sqrt{r}}. 
		\end{align*} 
	\end{lemma} 
	\begin{proof}
		See \cite[Lemma 4.5, page 72]{TI}. 
	\end{proof}

\begin{lemma}\label{beta}
For a natural number $k >1$, let $S = \{m^k : X \leq m^k \leq 2X\} \subset [X, 2X]$. Let $\beta_{\ell}(m)$ be an arithmetical function with
\begin{align*}
 |\beta_{\ell}(m)| =  \mathop{\sum_{m_1 \in S}\,\sum_{m_2 \in S}\,...\,\sum_{m_{2^{\ell+1}} \in S}}_{m \,=\,(m_1-m_2) + (m_3 - m_4) \,+\, ...\, + \,(m_{2^{\ell+1}-1}- m_{2^{\ell+1}}) } \,1\,. 
 \end{align*}
For $2^{\ell +1} > k \log k$, we have
\begin{align*}
||\beta_{\ell}||^2 =  \sum_{m}\,|\beta_{\ell}(m)|^2 \,\ll\,  X^{2^{\ell+2}/k - 1} .  
\end{align*}
\end{lemma}
\begin{proof}
We can write
\begin{align*}
 \sum_{m}\,|\beta_{\ell}(m)|^2 = \mathop{\mathop{\sum_{m}\,\sum_{m^k_1 \in S}\,\sum_{m^k_2 \in S}\,...\,\sum_{m^k_{2^{\ell+1}} \in S}\,1 \,\sum_{n^k_1 \in S}\,\sum_{n^k_2 \in S}\,...\,\sum_{n^k_{2^{\ell+1}} \in S} 
 1}_{m^k_1 - m^k_2 + m^k_3- m^k_4 + ... + m^k_{2^{\ell+1}-1} - m^k_{2^{\ell+1}} = m}}_{n^k_1 - n^k_2 + n^k_3 - n^k_4  ... + n^k_{2^{\ell+1}-1} - n^k_{2^{\ell+1}} = m}    
\end{align*}
\begin{align*}
    \ll&\,\,\, \mathop{\sum_{m}\,\sum_{m^k_1 \in S}\,\sum_{m^k_2 \in S}\,...\,\sum_{m^k_{2^{\ell+1}} \in S}\, \,\sum_{n^k_1 \in S}\,\sum_{n^k_2 \in S}\,...\,\sum_{n^k_{2^{\ell+1}} \in S}\, 
 1}_{m^k_1 + n^k_2 + m^k_3+ n^k_4 + ... + m^k_{2^{\ell+1}-1} + n^k_{2^{\ell+1}} = m = n^k_1 + m^k_2 + n^k_3 + m^k_4 + ... + n^k_{2^{\ell+1}-1} + m^k_{2^{\ell+1}}}
 =\,\,\, \mathop{\sum_{m}} |w_{\ell}(m)|^2,
\end{align*}
where $w_{\ell}(m)$ count the total number of representations of $m$ as the sum of $2^{\ell+1}$ many positive $k$-th powers. From \cite[Theorem $20.2$]{R1}, for $2^{\ell+1} > 2^k$, we have  
\begin{align*}
w_{\ell}(m) \asymp m^{2^{\ell+1}/k - 1} .   
\end{align*}
The condition $2^{\ell+1} > 2^k$ was further improved by Wooley to $2^{\ell+1} > k \log k$. Hence, we get
\begin{align*}
 \sum_{m}\,|\beta_{\ell}(m)|^2 \ll\, X^{2^{\ell+2}/k - 1}.   \end{align*}
 This proves our lemma.
\end{proof}
\vspace{0.3cm}

In the next Lemma, we will prove that, on average, the greatest common divisor factors will behave like $1$. 

\begin{lemma}\label{ll1}
Let $k \geq 1$ be a natural number, $X_1, X_2, ..., X_k \geq 1$  be any variables and $s_1, s_2, ..., s_k$ be any fixed $k$ non-zero integers. Let $f: \mathbb{N} \times \mathbb{N} \times ... \times \mathbb{N} \longrightarrow \mathbb{C}$ be any complex-valued function. For the real numbers $ 0 \leq r_1, r_2, ..., r_k \leq 1/2$, we have
\begin{align*}
    \sum_{n_1 \leq X_1}\,&\sum_{n_2 \leq X_2}\,...\sum_{n_k \leq X_k}\,(n_1, s_1)^{r_1}\,(n_2, s_2)^{r_2}\,...(n_k, s_k)^{r_k}\,|f(n_1, n_2, ..., n_k)|\\
    &\,\,\,\,\,\,\,\,\,\,\,\ll_{\epsilon, s_1, s_2,...,s_k}\,(X_1X_2...X_k)^{1/2 +\epsilon}\,\left(\sum_{n_1 \leq X_1}\,\sum_{n_2 \leq X_2}\,...\sum_{n_k \leq X_k}\,|f(n_1, n_2, ..., n_k)|^2\,\right)^{1/2}.
\end{align*}
\end{lemma}
\begin{proof}
Let \begin{align*}
    A(X_1, X_2, ..., X_k) =  \sum_{n_1 \leq X_1}\,\sum_{n_2 \leq X_2}\,...\sum_{n_k \leq X_k}\,(n_1, s_1)^{r_1}\,(n_2, s_2)^{r_2}&\,...(n_k, s_k)^{r_k}\,f(n_1, n_2, ..., n_k).
\end{align*}
On applying the Cauchy-Schwarz inequality to the sum over $n_1, n_2, ..., n_k$, we get
\begin{align*}
A(X_1, X_2, ..., X_k)\,\ll\, S_1 \,S_2\,...S_k\,\left(\sum_{n_1 \leq X_1}\,\sum_{n_2 \leq X_2}\,...\sum_{n_k \leq X_k}\,|f(n_1, n_2, ..., n_k)|^2 \right)^{1/2},    
\end{align*}
where 
\begin{align*}
    S^2_i = \sum_{n_i \leq X_i}\, |(n_i, s_i)^{r_i}|^2,\,\,\,\,\,\,\,\,\text{for all}\,\,\, 1 \leq i \leq k.
\end{align*}
Let $(n_i, s_i) = \ell_i$, then we have
\begin{align*}
\sum_{n_i \leq X_i}\,|(n_i, s_i)^{r_i}|^2 \, \leq \,\sum_{\ell_i | s_i}\,\,\sum_{n_i \leq X_i/\ell_i}\,\ell^{2r_i}_i \,\leq\, X_i\,\sum_{\ell_i | s_i} \, \ell^{2r_i -1}_i\,\ll_{\epsilon_i}\, {(s^{2r_i -1}_i X_i)^{1+\epsilon_i}} 
\end{align*}
As $2r_i - 1 \leq 1$, so we have
\begin{align*}
 S_i\,\ll_{\epsilon_i, s_i}\, \left(X^{1 + \epsilon_i}_i\right)^{1/2},   
\end{align*}
where $\epsilon_i > 0$'s are arbitrarily small. Thus, on plugging the bound for $S_i$, we get
\begin{align*}
A(X_1, X_2, ..., X_k) \,&\ll_{\epsilon, s_1, s_2,...,s_k}\,\\
&\,(X_1X_2...X_k)^{1/2 +\epsilon}\,\left(\sum_{n_1 \leq X_1}\,\sum_{n_2 \leq X_2}\,...\sum_{n_k \leq X_k}\,|f(n_1, n_2, ..., n_k)|^2\,\right)^{1/2}.
\end{align*}
This proves our lemma.
\end{proof}

\section{Proof of Theorem $1$}

Let $\alpha(n)$ and $\beta(n)$ denote two sequences of complex numbers, and $A_{\pi}(1,n)$ represent the normalized $n$-th Fourier coefficients of the $SL(3,\mathbb{Z})$ Hecke-Maass cusp form $\pi$. We want to get the cancellations in the following sum:

\begin{align*}
    S(N,M)  =\, \sum_{1 \leq n \leq N}\,\sum_{1 \leq m \leq M}\, \alpha(n)\,\beta(m)\, A_{\pi}(1,n+m),
\end{align*}
where $M < N^{1-\epsilon}$. On applying the Cauchy-Schwarz inequality to the sum over $n$, we get that $S(N, M)$ is dominated by
\begin{align*}
		||\alpha||_2 \, \left(\sum_{1 \leq n \leq N}\,W_1\left(\frac{n}{N}\right)\left|\sum_{1 \leq m \leq M} \beta(m)\,A_{\pi}(1,n+m)\right|^2\right)^{1/2}  =  ||\alpha||_2 \, \left(\gamma(N, M)\right)^{1/2} ,
	\end{align*}
		where $\gamma(n,m)$ is given by 
\begin{align*}
		\gamma(N, M) = \sum_{1 \leq m_1 \leq M} \,\sum_{1 \leq m_2 \leq M}\, \beta(m_1)\,\overline{\beta(m_2)} \, \sum_{n=1}^\infty\,W_1\left(\frac{n}{N   }\right)\, A_{\pi}(1,n+m_1)\,\overline{A_{\pi}(1,n+m_2)},     
\end{align*}
and $W_1$ is a smooth bump function supported on $[-2,2]$, that takes positive values with $W_1(x) = 1$ on $[-1,1]$, satisfying $x^j W_1^{(j)}(x) \ll 1$ for all $j \geq 1$. Using the change of variables $n+m_1 \longrightarrow n$ and $m_2 - m_1 = m$, we can write $\gamma(N, M) $ as
	\begin{align*}
		\sum_{1 \leq m_1 \leq M}\, \sum_{(1-m_1) \leq m \leq (M-m_1)} \, \beta(m_1)\,\overline{\beta(m+m_1)} \, \sum_{n=1}^\infty\,W_1\left({\frac{n-m_1}{N}}\right) A_{\pi}(1,n)\,\overline{A_{\pi}(1,n+m)},    
	\end{align*}
Now we apply the Cauchy-Schwarz inequality to the sum over $m$, and then $m_1$, we arrive at  
\begin{align*}
	\gamma(N, M)\,&\ll\,	\,||\beta||^2_2\,\sqrt{M}  \left(\sum_{m \in \mathbb{Z}}\,W\left(\frac{m}{M   }\right)\,\left|\sum_{n=1}^\infty\,W_1\left(\frac{n}{N_1}\right) A_{\pi}(1,n) \overline{A_{\pi}(1,n+m)}\right|^2 \right)^{1/2}\\
 \,&\ll\,||\beta||^2_2 \sqrt{M} (\eta(N, M)^{1/2},
\end{align*}
 where $W$ is a non-negative smooth even function supported on $[-2,2]$ with $W(x)=1$ on $[-1,1]$, satisfying $x^j W^{(j)}(x) \ll 1$ and $N_1 :=  M+N \asymp  N$. Also, we have
	\begin{align*}
	\eta(N, M)=& \,\sum_{n_1 =1}^\infty\,\sum_{n_2 =1}^\infty W_1\left({\frac{n_1}{N_1}}\right) W_2\left({\frac{n_2}{N_1}}\right) A_{\pi}(1,n_1) \overline{A_{\pi}(1,n_2)}\\
    &\,\times\,\sum_{m \in \mathbb{Z}} W\left(\frac{m}{M   }\right) A_{\pi}(1, n_1 +m) \,\overline{A_{\pi}(1, n_2 +m)}, 
	\end{align*}
  where $W_2$ is another smooth bump function satisfying the same properties as $W_1$, which is introduced to symmetrize the $n_1, n_2$ variables after Cauchy–Schwarz. Again, using the change of the variable $n_1 + m \longrightarrow n_3$, we can write 
\begin{align*}
\eta(N, M) &=  \sum_{n_1=1}^\infty\,\sum_{n_2=1}^\infty\,W_1\left({\frac{n_1}{N_1}}\right) W_2\left({\frac{n_2}{N_1}}\right)\, A_{\pi}(1,n_1)\,\overline{A_{\pi}(1,n_2)}\,\\
&\times\, \sum_{n_3 \in \mathbb{Z}}\,W\left(\frac{n_3-n_1}{M}\right)\, A_{\pi}(1, n_3)\, \overline{A_{\pi}(1, n_2 +n_3-n_1)}.
\end{align*}
Here we note that after the change of variables $n_3 = n_1 + m$, the support conditions on the
weights $W_1(n_1/N_1)$, $W(m/M)$, together with ${A}_{\pi}(1, n_3) = 0$ for $n_3 \leq 0$ imply that $1 \le n_3 \le 2N_1 + 2M$.
We introduce a smooth weight function
$W_3 \in C_c^\infty(\mathbb{R})$ such that $W_3(x)=1$ for $x \in [-1,1]$,
$\operatorname{supp}(W_3) \subset [-2,2]$, satisfying $x^j W_3^{(j)}(x) \ll 1$, and define
$
N_2 := 2N_1 + 2M \asymp N$. Now we note that $W_3(n_3/N_2) \equiv 1$ on the effective support of the summation, so the insertion of this weight $W_3$ does not change the value of the sum. Therefore, we can write
\begin{align*}
\eta(N, M) &=  \sum_{n_1=1}^\infty\,\sum_{n_2=1}^\infty\,W_1\left({\frac{n_1}{N_1}}\right) W_2\left({\frac{n_2}{N_1}}\right)\, A_{\pi}(1,n_1)\,\overline{A_{\pi}(1,n_2)}\,\\
&\times\, \sum_{n_3 \in \mathbb{Z}}\, W\left(\frac{n_3-n_1}{M}\right)\,W_3\left({\frac{n_3}{N_2}}\right) \,A_{\pi}(1, n_3)\,\overline{A_{\pi}(1, n_2 +n_3-n_1)},
\end{align*}
 The factor $W\left({(n_3-n_1)}/{M}\right)$ records the information on the change of variables and aids in obtaining additional savings later. Now, we detect the equation $n_2+n_3-n_1 = r$ with the help of the delta symbol and using the support conditions of $W_1, W_2$ and $W_3$, we introduce the weight function $W_4$ to the $n_4$-sum satisfying the similar properties and support conditions as $V_3$. With all this, we can write
\begin{align*}
\eta(N, M) \,=	\,& \,\sum_{n_1 = 1}^{\infty}\,A_{\pi}(1,n_1)\,W_1\left(\frac{n_1}{N_1}\right)\,\sum_{n_2 = 1}^{\infty}\,\overline{A_{\pi}(1,n_2)}\,W_2\left(\frac{n_2}{N_1}\right)\\
  &\times \, \sum_{n_3 =1}^{\infty}\, A_{\pi}(1, n_3)\,W\left(\frac{n_3-n_1}{M}\right)\,W_3\left(\frac{n_3}{N_2}\right)\\
  &\times \, \sum_{r = 1}^{\infty}\,\overline{A_{\pi}(1, r)}\,W_4\left(\frac{r}{N_3}\right) \,\delta({r= n_2 + n_3 -n_1 }),
\end{align*}
\vspace{0.1cm}

\noindent where $N_3 =N_2$. Using the expansion of the delta symbol $\delta(n, m)$ due to Duke, Friedlander, and Iwaniec from the section \ref{a1}, with $q \sim Q \asymp \sqrt{N}$, the expression $\eta(n,m)$ can be written as
	\begin{align}\label{L13}
		\eta(N, M)  =&\, \frac{1}{Q}\,\sum_{q\leq Q}\frac{1}{q}\, \sideset{}{^\star} \sum_{a\bmod q}\,\int_{\mathbb{R}} \, \psi(q,u)\notag\\
    &\times\, \sum_{n_1 =1}^{\infty}\,A_{\pi}(1,n_1) \,e\left(\frac{an_1}{q}\right)\,e\left(\frac{n_1u}{qQ}\right)\,W_1\left(\frac{n_1}{N_1}\right) \notag\\ 
    &\times\, \sum_{n_2 =1}^{\infty}\,\overline{A_{\pi}(1,n_2)}\, e\left(\frac{-an_2}{q}\right)\,e\left(\frac{-n_2u}{qQ}\right)\,W_2\left(\frac{n_2}{N_1}\right) \notag\\
    &\times\, \sum_{n_3 =1}^{\infty}\,{A_{\pi}(1,n_3)}\, e\left(\frac{-an_3}{q}\right)\,e\left(\frac{-n_3u}{qQ}\right)\,W_3\left(\frac{n_3}{N_2}\right)\,W\left(\frac{n_3-n_1}{M}\right) \notag\\
    &\times\, \sum_{r =1}^{\infty}\,\overline{A_{\pi}(1,r)}\, e\left(\frac{ar}{q}\right)\,e\left(\frac{ru}{qQ}\right)\,W_4\left(\frac{r}{N_3}\right)du + O(N^{-2025}). 
	\end{align}
Here, $\psi$ is a bump function supported on $[-2 N^{\epsilon}, 2 N^{\epsilon}]$ such that $\psi(u) = 1$, for all $u \in [- N^{\epsilon},  N^{\epsilon}]$, with $\psi^{(j)}(u) \ll_j 1$. We will now apply the GL$(3)$ Voronoi summation formula to the $n_1, n_2, n_3$, and $r$ sums, starting with the $n_3$-sum. The following lemmas provide the necessary details.

\begin{lemma}\label{L14}
Let $\mathcal{S}_3$ denote the $n_3$-sum as given in \eqref{L13}. Then, we have
\begin{align*}
    \mathcal{S}_3 =\, \frac{N^{2/3}_2}{q}\, \sum_{\pm}\, \sum_{n'_{3}|q}(n'_{3})^{1/3} \,\sum_{(n'_3)^2n''_{3}\ll M_0}\,  \frac{A_{\pi}(n''_{3},n'_{3})}{(n''_{3})^{1/3}} \,S\left( -\bar{a}, \pm n''_{3}; q/n'_{3}\right)\, \mathcal{I}^{\pm}_3 + O(N^{-A}),
\end{align*}
where $ M_0= \text{max}\left\{{q^3}/{N_2}, {(qN_2)^3}/{M^3 N_2}, N^{1/2}_2 |u|^3\right\}N^{\epsilon}$, $A$ is any positive real number, and
\begin{align*}
    \mathcal{I}^{\pm}_3 =\, \int_{0}^\infty\, W_3(x_3)\,W\left(\frac{N_2x_3-n_1}{M}\right)\, e\left(-\frac{N_2ux_3}{qQ} \pm  \frac{3(N_2x_3 (n'_3)^2n''_{3})^{1/3}}{q} \right)\,dx_3 .
\end{align*}
Here $W_3$ is a new weight function different from the one as in equation \eqref{L13}.
\end{lemma}
\begin{proof}
We have,
\begin{align*}
    \mathcal{S}_3 =\, \sum_{n_3 =1}^{\infty}\, A_{\pi}(1,n_3)\, e\left(\frac{-an_3}{q}\right)\,e\left(\frac{-n_3u}{qQ}\right)\,W_3\left(\frac{n_3}{N_2}\right)\,W\left(\frac{n_3-n_1}{M}\right).
\end{align*}
Applying the GL$(3)$-Voronoi summation formula from Lemma \ref{gl3voronoi} to the sum $\mathcal{S}_3$, we get
\begin{align} \label{A10}
    \mathcal{S}_3 =\, q\, \sum_{\pm}\, \sum_{n'_{3}|q} \,\sum_{n''_{3}=1}^{\infty}\, \frac{A_{\pi}(n''_{3},n'_{3})}{n'_{3} n''_{3}}\, S\left( -\bar{a}, \pm n''_{3}; q/n'_{3}\right)\, \mathcal{W}^{\pm}_3 \left(\frac{(n'_{3})^2 n''_{3}}{q^3 }\right).
\end{align} 
Using the asymptotic expansion of the integral transform $\mathcal{W}^{\pm}_3 \left(\frac{(n'_{3})^2 n''_{3}}{q^3 }\right)$ from Lemma \ref{GL3oscilation} and simplifying it, we can write
\begin{align*}
    \mathcal{W}^{\pm}_3 \left(\frac{(n'_{3})^2 n''_{3}}{q^3 }\right) =\, \frac{N^{2/3}_2}{q^2}\, ((n'_{3})^2 n''_{3})^{2/3}\, \mathcal{I}^{\pm}_3 + \text{lower order terms},
\end{align*} 
where
\begin{align*}
    \mathcal{I}^{\pm}_3 =\, \int_{0}^\infty\, W_3(x_3)\,W\left(\frac{N_2x_3-n_1}{M}\right)\, e\left(-\frac{N_2ux_3}{qQ} \right)\, e\left(\pm  \frac{3(N_2x_3 (n'_3)^2n''_{3})^{1/3}}{q} \right)\,dx_3 .  
\end{align*}
Here, $W_3$ is a new weight function different from the one in the above expression of $\mathcal{S}_3$. For the abuse of notation, we are denoting it with the same notation. Now, using integration by parts $j$-times, we arrived at the following expression 
\begin{align*}
    \mathcal{I}^{\pm}_3 \ll \left(1 + \frac{N_2}{M} + \frac{N_2|u|}{qQ}\right)^j \left(\frac{q}{(N_2(n'_3)^2n''_3)^{1/3}}\right)^j.
\end{align*}
Hence, the integral transform $\mathcal{I}^{\pm}_3$ will be negligibly small unless 
\begin{align*}
    (n'_3)^2n''_3 \ll \text{max}\left\{\frac{q^3}{N_2}, \frac{(qN_2)^3}{M^3 N_2}, N^{1/2}_2 |u|^3\right\} =: M_0.
\end{align*}
Thus, the integral transform $\mathcal{W}^{\pm}_3 \left(\frac{(n'_{3})^2 n''_{3}}{q^3 }\right)$ is negligibly small if $(n'_3)^2n''_3 \gg {M_0}$. Ignoring the lower-order terms (see remark \ref{rm24}) and using equation \eqref{A10}, we can write 
\begin{align}\label{A16}
    \mathcal{S}_3 = \frac{N^{2/3}_2}{q} \sum_{\pm} \sum_{n'_{3}|q} \sum_{(n'_3)^2n''_{3}\ll M_0}  \frac{A_{\pi}(n''_{3},n'_{3})}{(n'_{3})^{-1/3}(n''_{3})^{1/3}} S\left( -\bar{a}, \pm n''_{3}; q/n'_{3}\right) \mathcal{I}^{\pm}_3 + O(N^{-A}).
\end{align}
This proves our lemma.
\end{proof}

\noindent Next, we will apply the $GL(3)$ Voronoi summation formula to the sum over $n_1$. We have the following lemma.
\vspace{0.3cm}

\begin{lemma}\label{L15}
Let $\mathcal{S}_1$ denote the $n_1$-sum be as given in equation \eqref{L13} i.e.
  \begin{align*}
\mathcal{S}_1\,=\,  \sum_{n_1 =1}^{\infty}\,A_{\pi}(1,n_1)\,   \,e\left(\frac{an_1}{q}\right)\,e\left(\frac{n_1u}{qQ}\right)\,W_1\left(\frac{n_1}{N_1}\right)\,W\left(\frac{N_2x_3-n_1}{M}\right).  
\end{align*}
Then we have 
\begin{align*}
		\mathcal{S}_1 = \frac{N^{2/3}_1}{q} \sum_{\pm} \sum_{n'_{1}|q}\,(n'_{1})^{1/3} \sum_{(n'_1)^2n''_{1}\ll M_0}  \frac{A_{\pi}(n''_{1},n'_{1})}{(n''_{1})^{1/3}} S\left(\bar{a}, \pm n''_{1}; q/n'_{1}\right)\,\, \mathcal{I}^{\pm}_1 + O(N^{-A}),
	\end{align*} where $ M_0= \text{max}\left\{{q^3}/{N_1}, {(qN_1)^3}/{M^3 N_1}, N^{1/2}_1 |u|^3\right\}N^{\epsilon}$, $A$ is any positive real number and,
	\begin{align*}
		\mathcal{I}^{\pm}_1\, = \,\int_{0}^\infty W_1(x_1)\,W\left(\frac{N_2x_3- N_1 x_1}{M}\right)\, e\left(\frac{N_1 ux_1}{qQ} \pm  \frac{3(N_1 x_1 (n'_1)^2n''_{1})^{1/3}}{q} \right)\,dx_1 .
	\end{align*}
 Here $W_1$ is a new weight function different from the one as in equation \eqref{L13}.
 \end{lemma}
\begin{proof}
Proof will follow the same steps as in Lemma \ref{L14}.    
\end{proof}
\vspace{0.3cm}

 \begin{lemma}\label{L16}
Let $\mathcal{S}_2$ denote the $n_2$-sum be as given in equation \eqref{L13}.
Then we have 
\begin{align*}
		\mathcal{S}_2 = \frac{N^{2/3}_1}{q} \sum_{\pm} \sum_{n'_{2}|q}\,(n'_{2})^{1/3} \sum_{(n'_2)^2n''_{2}\ll N_0}  \frac{\overline{A_{\pi}(n''_{2},n'_{2})}}{(n''_{2})^{1/3}} S\left(-\bar{a}, \pm n''_{2}; q/n'_{2}\right)\,\, \mathcal{I}^{\pm}_2 + O(N^{-A}),
	\end{align*} where $ N_0= \text{max}\left\{{q^3}/{N_1}, N^{1/2}_1 |u|^3\right\}N^{\epsilon}$, $A$ is any positive real number and,
	\begin{align*}
		\mathcal{I}^{\pm}_2\, = \,\int_{0}^\infty W_2(x_2)\, e\left(-\frac{N_1ux_2}{qQ} \pm  \frac{3(N_1 x_2 (n'_2)^2n''_{2})^{1/3}}{q} \right)\,dx_2 .
	\end{align*}
 Here $W_2$ is a new weight function different from the one as in equation \eqref{L13}.
 \end{lemma}
\begin{proof}
Proof will follow the same steps as in Lemma \ref{L14}.    
\end{proof}
\vspace{0.3cm}

 \begin{lemma}\label{L17}
Let $\mathcal{S}_4$ denote the $r$-sum be as given in equation \eqref{L13}.
Then we have 
\begin{align*}
		\mathcal{S}_4 = \frac{N^{2/3}_3}{q} \sum_{\pm} \sum_{r_1|q}\,r^{1/3}_1 \sum_{r^2_1 r_2\ll N_0}  \frac{\overline{A_{\pi}(r_2,r_1)}}{r^{1/3}_2} S\left(-\bar{a}, \pm r_{2}; q/r_1\right)\,\, \mathcal{I}^{\pm}_4 + O(N^{-A}),
	\end{align*} where $ N_0= \text{max}\left\{{q^3}/{N_3}, N^{1/2}_3 |u|^3\right\}N^{\epsilon}$, $A$ is any positive real number and,
	\begin{align*}
		\mathcal{I}^{\pm}_4\, = \,\int_{0}^\infty W_4(x_4)\, e\left(\frac{N_3 ux_4}{qQ} \pm  \frac{3(N_3r^2_1r_2)^{1/3}}{q} \right)\,dx_4.
	\end{align*}
 Here $W_4$ is a new weight function different from the one as in equation \eqref{L13}.
 \end{lemma}
\begin{proof}
Proof will follow the same steps as in Lemma \ref{L14}.    
\end{proof}

\noindent Using Lemma \ref{L14}, Lemma \ref{L15}, Lemma \ref{L16}, and Lemma \ref{L17} to the equation \eqref{L13}, $\eta(N, M)$ will get transformed into the following expression
\begin{align}\label{s1}
\eta(N, M) = \,\,&\frac{N^{8/3}}{Q}\sum_{q\leq Q}\frac{1}{q^{5}}\,\int_{|u| \leq Q^{\epsilon}} \, \psi(q,u) \,\,  \sum_{\pm}\,\sum_{n'_{1}|q}\,(n'_{1})^{1/3}\, \sum_{(n'_1)^2n''_{1}\ll M'_0}\,  \frac{{A_{\pi}(n''_{1},n'_{1})}}{ (n''_{1})^{1/3}} \notag\\
 &\times \sum_{\pm}\,\sum_{n'_{2}|q}\,(n'_{2})^{1/3}\, \sum_{(n'_2)^2n''_{2} \ll N'_0}\,  \frac{\overline{A_{\pi}(n''_{2},n'_{2})}}{ (n''_{2})^{1/3}}    \,\sum_{\pm}\,\sum_{n'_{3}|q}\,(n'_{3})^{1/3}\, \sum_{(n'_3)^2n''_{3} \ll M'_0}\,  \frac{A_{\pi}(n''_{3},n'_{3})}{ (n''_{3})^{1/3}}\, \notag \\
 &\times \,\,\sum_{\pm} \sum_{r_{1}|q}\,r^{1/3}_{1} \, \sum_{r^2_1r_{2}\ll N'_0}  \frac{\overline{A_{\pi}(r_{2},r_{1})}}{ r^{1/3}_{2}}\,\, \mathcal{C}(...)\, \,\mathcal{I}(...)\,du + O(X^{-A}),
\end{align}
where
\begin{align*}
  &M_0 \asymp M'_0 = {\max}\left\{q^3/N, (qN)^3/M^3N, N^{1/2}|u|^3\right\}N^{\epsilon}\\
  & N_0 \asymp N'_0 = {\max}\left\{q^3/N, N^{1/2}|u|^3\right\}N^{\epsilon},
\end{align*}
the character sum $\mathcal{C}(...)$ is given by
\begin{align}\label{L10}
\mathcal{C}(...) = \sideset{}{^\star} \sum_{a \bmod q}\, S\left(\overline{a}, \pm n''_{1}; q/n'_{1}\right) S\left(-\overline{a}, \pm n''_{2}; q/n'_{2}\right) S\left(-\overline{a}, \pm n''_{3}; q/n_{3}\right) S\left(\overline{a}, \pm r_{2}; q/r_{1}\right),  
\end{align}
and the integral transform $\mathcal{I}(...)$ is given by
\begin{align}\label{L11}
\mathcal{I}(...) = &\,\int_{0}^\infty W_1(x_1)\, e\left(\frac{N_1ux_1}{qQ} \pm  \frac{3(N_1(n'_1)^2 n''_1 x_1)^{1/3}}{q} \right)\notag\\
&\times \,\int_{0}^\infty W_2(x_2)\, e\left(-\frac{N_1ux_2}{qQ} \pm  \frac{3(N_1(n'_2)^2 n''_2 x_2)^{1/3}}{q} \right)\,\notag\\
&\times \,\int_{0}^\infty W_3(x_3)\,W\left(\frac{N_2x_3- N_1x_1}{M}\right) \,\, e\left(- \frac{N_2ux_3}{qQ} \pm  \frac{3(N_2 (n'_3)^2 n''_3 x_3)^{1/3}}{q} \right)\notag\\
&\times\,\int_{0}^\infty W_4(x_4)\, e\left( \frac{N_3ux_4}{qQ} \pm  \frac{3(N_3 r_1^2 r_2 x_4)^{1/3}}{q} \right)dx_1\,dx_2\,dx_3\,dx_4.
\end{align}
\vspace{0.5cm}

\subsection{Analysis of the character sum}\label{M26}
In this section, our objective is to determine an upper bound for the character sum
\begin{align*}
\mathcal{C}(...) = \sideset{}{^\star} \sum_{a \bmod q}\, S\left(\overline{a}, \pm n''_{1}; q/n'_{1}\right) S\left(-\overline{a}, \pm n''_{2}; q/n'_{2}\right) S\left(-\overline{a}, \pm n''_{3}; q/n'_{3}\right) S\left(\overline{a}, \pm r_{2}; q/r_{1}\right).  
\end{align*}
Let $q = q_1 q_2 q_3$, where $q_1$ is the largest factor of $q$ such that $q_1|[n'_1, n'_2, n'_3, r_1]$ with $(q_1, q/q_1) = 1$,  $q_2$ is the largest factor of $q/q_1$ such that $q_2|(n'_1 n'_2 n'_3 r_1)^{\infty}$ with $(q_2, q/q_1q_2) =1$ and $q_3 = q/q_1q_2$. Here $[s, t]$ means the least common multiple of $s$ and $t$. Notice that $q_2$ is square full and  $(n'_1 n'_2 n'_3 r_1)| q_1q_2$. Then we can write any reduced residue class $a$ modulo $q$  as $ a = a_1q_3 \overline{q_3} + a_2 q_1q_2\,\overline{q_1q_2}$, where $a_1$ and $a_2$ are the reduced residue classes modulo $q_1q_2$ and $q_3$, respectively. Also, we have 
\begin{align*}
 q_3\overline{q_3} \equiv
 1 \,\mathrm{mod} \, q_1q_2, \hspace{1cm}  (q_1q_2)\overline{(q_1q_2)} \equiv
 1 \,\mathrm{mod} \, q_3.  
\end{align*}
  By changing the variable $\overline{a} \longrightarrow a$ and using the following well-known properties of Kloosterman sums
\begin{align}\label{B31}
		S (a, b; c_1 c_2 ) = S (a \overline{c_2}, b  \overline{c_2}; c_1  ) S (a \overline{c_1}, b  \overline{c_1}; c_2  )   \ \ \ \ \ \ \ \ \ \ \ \  \textrm{and,}
	\end{align}	
	\begin{align}\label{B32}
		S (am, b m; c)  =  S (a m^2, b ; c) =  S (a, b m^2; c)  \ \ \  \ \ \ \textrm{if} \ \ (m, c) = 1,
	\end{align}
we can decompose our character sum into two parts as follows
\begin{align}\label{m4}
 |\mathcal{C}(...)| = \,\mathcal{C}_1(...)\,\mathcal{C}_2(...),  
\end{align}
where
\begin{align*}
 \mathcal{C}_1(...) = &\sideset{}{^\star} \sum_{a_1 \bmod q_1q_2}\, S\left(a_1, \pm n''_{1}\overline{q_3}^2; q_1q_2/n'_{1}\right) S\left(-a_1, \pm n''_{2}\overline{q_3}^2; q_1q_2/n'_{2}\right)\\
&\,\,\,\,\,\,\, \times\,S\left(-a_1, \pm n''_{3}\overline{q_3}^2; q_1q_2/n'_{3}\right) S\left(a_1, \pm r_{2}\overline{q_3}^2; q_1q_2/r_{1}\right)
 \end{align*}
 and
 \begin{align*}
  \mathcal{C}_2(...) = &\,\sideset{}{^\star} \sum_{a_2 \bmod q_3}\, S\left(a_2, \pm n''_{1}(\overline{q_1q_2/n'_1})^2; q_3\right) S\left(-a_2, \pm n''_{2}(\overline{q_1q_2/n'_2})^2; q_3\right)\\
  &\,\,\,\,\,\,\,\,\,\times\,S\left(-a_2, \pm n''_{3}(\overline{q_1q_2/n'_3})^2; q_{3}\right) S\left(a_2, \pm r_{2}(\overline{q_1q_2/r_1})^2; q_3\right).  
\end{align*}
Estimating the $a_1$-sum trivially and using the Weil's bound for the Kloosterman sums in $\mathcal{C}_1(...) $, we obtain
\begin{align*}
\mathcal{C}_1(...) \,& \ll \, \frac{(q_1q_2)^3}{\sqrt{n'_1 n'_2n'_3r_1}}\,d\left(\frac{q_1q_2}{n'_1}\right)\,d\left(\frac{q_1q_2}{n'_2}\right)\,d\left(\frac{q_1q_2}{n'_3}\right)\,d\left(\frac{q_1q_2}{r_1}\right)\,\left(a_1,  n''_{1}\overline{q_3}^2, \frac{q_1q_2}{n'_1}\right)^{1/2}\\
&\,\,\,\,\,\,\,\times \,\left(a_1,  n''_{2}\overline{q_3}^2, \frac{q_1q_2}{n'_2}\right)^{1/2}\,\left(a_1,  n''_{3}\overline{q_3}^2, \frac{q_1q_2}{n'_3}\right)^{1/2}\,\left(a_1,  r_{2}\overline{q_3}^2, \frac{q_1q_2}{r_1}\right)^{1/2}.
\end{align*}
Note that, as $(a_1, q_1q_2/n_1) = 1 = (\overline{q_3}^2, q_1q_2/n_1)$. So, we have
\begin{align*}
\left(a_1,  n''_{1}\overline{q_3}^2, \frac{q_1q_2}{n'_1}\right) = \left(n''_{1}\overline{q_3}^2, \frac{q_1q_2}{n'_1}\right) =  \left(n''_{1}, \frac{q_1q_2}{n'_1}\right)   
\end{align*}
Using above calculations and the well known estimates $d(n) \ll_{\epsilon} n^{\epsilon}$, we can write
\begin{align*}
\mathcal{C}_1(...) \,& \ll_{\epsilon} \, \frac{(q_1q_2)^3}{\sqrt{n'_1 n'_2n'_3r_1}}\,\,\left(n''_{1}, \frac{q_1q_2}{n'_1}\right)^{1/2}\,\left(n''_{2}, \frac{q_1q_2}{n'_2}\right)^{1/2}\,\left(n''_{3}, \frac{q_1q_2}{n'_3}\right)^{1/2}\,\left(r_{2}, \frac{q_1q_2}{r_1}\right)^{1/2} \left(\frac{q_1q_2}{n_1}\right)^{\epsilon}.    
\end{align*}
In the Lemma $\ref{ll1}$, we have proved that when we consider the greatest common divisor factors on average, they will not contribute anything to the final bound other than $X^{\epsilon}$ and behave like $1$. So, when we use the above bound for $\mathcal{C}_1(...)$ in equation \eqref{s1} and consider these greatest common divisor factors on average (i.e., with the $n''_1, n''_2, n''_3$ and $r_2$ summations), we can ignore these factors. So, to avoid writing a messy bound, let us write
\begin{align}\label{m5}
\mathcal{C}_1(...) \,\ll_{\epsilon} \, \frac{(q_1q_2)^3}{\sqrt{n'_1 n'_2n'_3r_1}}\,(q_1q_2)^{\epsilon}.
\end{align}

\noindent To analyze the character sum $\mathcal{C}_2(...)$ and find its bound, we further break $q_3$ to extract its square full part. We write $q_3 = q'_3q''_3$ such that $(q'_3, 2q''_3) =1$, where $q'_3$ is square free and $4q''_3$ is square full. In this situation, we can rewrite the character sum $\mathcal{C}_2(...)$ as
\begin{align*}
 \mathcal{C}_2(...) = \,\,\mathcal{C}'_2(...)\,\mathcal{C}''_2(...),  
\end{align*}
where
\begin{align*}
\mathcal{C}'_2(...) =  &\,\sideset{}{^\star} \sum_{a'_2\,\mathrm mod\, q'_3}\, S\left(a'_2, \pm n''_{1}(\overline{q_1q_2q''_3/n'_1})^2; q'_3\right) S\left(-a'_2, \pm n''_{2}(\overline{q_1q_2q''_3/n'_2})^2; q'_3\right)\\
  &\,\,\,\,\,\,\,\,\,\times\,S\left(-a'_2, \pm n''_{3}(\overline{q_1q_2q''_3/n'_3})^2; q'_{3}\right) S\left(a'_2, \pm r_{2}(\overline{q_1q_2q''_3/r_1})^2; q'_3\right),     
 \end{align*}
\begin{align*}
\mathcal{C}''_2(...) =  &\,\sideset{}{^\star} \sum_{a''_2\,\mathrm mod\, q''_3}\, S\left(a''_2, \pm n''_{1}(\overline{q_1q_2q'_3/n'_1})^2; q''_3\right) S\left(-a''_2, \pm n''_{2}(\overline{q_1q_2q'_3/n'_2})^2; q''_3\right)\\
  &\,\,\,\,\,\,\,\,\,\times\,S\left(-a''_2, \pm n''_{3}(\overline{q_1q_2q'_3/n'_3})^2; q''_{3}\right) S\left(a''_2, \pm r_{2}(\overline{q_1q_2q'_3/r_1})^2; q''_3\right).     
 \end{align*}
Again, estimating $\mathcal{C}''_2(...)$ using Weil's bound for Kloosterman sums with Lemma \ref{ll1} as for \eqref{m5}, we get
\begin{align*}
\mathcal{C}''_2(...) \,\ll_{\epsilon}\, (q''_3)^{3+\epsilon}.     
\end{align*}
We now proceed to estimate $\mathcal{C}'_2(...) $. Since we have chosen $q'_3$ as square free such that its prime factorization looks like $q'_3 = p_1p_2....p_k,\, k \geq 1$. So, we can rewrite $\mathcal{C}'_2(...)$ as
\begin{align}\label{m2}
 \mathcal{C}'_2(...) = \prod_{1 \leq i \leq k}\, \mathcal{C}'_{2,i}(...), \hspace{1cm} \text{where}   
\end{align}
\begin{align*}
\mathcal{C}'_{2,i}(...) =   &\,\sideset{}{^\star} \sum_{x\,\mathrm mod\, p_i}\, S\left(x, \pm n''_{1}(\overline{q_1q_2q'_3p'_i/n'_1})^2; p_i\right) S\left(-x, \pm n''_{2}(\overline{q_1q_2q'_3p'_i/n'_2})^2; p_i\right)\\
  &\,\,\,\,\,\,\,\,\,\times\,S\left(-x, \pm n''_{3}(\overline{q_1q_2q'_3p'_i/n'_3})^2; p_i\right) S\left(x, \pm r_{2}(\overline{q_1q_2q'_3p'_i/r_1})^2; p_i\right).      
\end{align*}
Here $p'_i = q'_3/p_i$. Define $s_j = (\overline{q_1q_2q'_3p'_i/n'_j})^2,\, 1 \leq  j \leq 3$ and $s_4 = (\overline{q_1q_2q'_3p'_i/r_1})^2$, where $(s_j, p_i) = 1$ for all $1 \leq j \leq 4,\, \, 1 \leq i \leq k$. We have
\begin{align}\label{m1}
\mathcal{C}'_{2,i}(...) \,= \,\sideset{}{^\star} \sum_{x\,\mathrm mod\, p_i}\, S\left(x, n''_{1}s_1; p_i\right)\,S\left(x, n''_{2}s_2; p_i\right)\,S\left(x, n''_{3}s_3; p_i\right)\,S\left(x, r_2s_4; p_i\right)   
\end{align}
Here we are replacing $\pm$ to $+$ in every Kloosterman sum to analyze only one such term out of $16$-many as the analysis will be similar for the other remaining terms. We will estimate the above character sum $\mathcal{C}'_{2, i}(...) $ in the next lemma.
\vspace{0.3cm}

\begin{lemma}\label{t1}
We have
\begin{align*}
 \mathcal{C}'_{2,i}(...) \,\ll\, p^{5/2}_i.    
\end{align*}
\end{lemma}
\begin{proof}
We can rephrase the character sum presented in equation \eqref{m1} as
\begin{align*}
\mathcal{C}'_{2,i}(...) \,=\, \sideset{}{^\star} \sum_{x\,\mathrm{mod}\, p_i}\, S\left(1, \gamma_1(x); p_i\right)\,S\left(1, \gamma_2(x); p_i\right)\,S\left(1, \gamma_3(x)x; p_i\right)\,S\left(1, \gamma_4(x); p_i\right),  
\end{align*}
where
\begin{align*}
    \gamma_i = \begin{pmatrix}
        n''_is_i& \,0\\
        0& \, 1
    \end{pmatrix}\,\,\, \text{for}\,\, i = 1,2,3 \,\,\,\,\, \text{and}\,\,\,\, \gamma_4 = \begin{pmatrix}
        r_2s_4& \,0\\
        0& \, 1
    \end{pmatrix}.
\end{align*}
Also, 
\begin{align*}
    \text{for} \,\,\,\gamma = \begin{pmatrix}
        a& \,b\\
        c& \, d
    \end{pmatrix} \in {PGL_2}(\mathbb{F}_p), \,\,\text{we mean}\,\,\,\, \gamma(x)  = \begin{pmatrix}
        a& \,b\\
        c& \, d
    \end{pmatrix} (x) \,=\, \frac{ax+b}{cx+d}.
\end{align*}
Fouvry, Kowalski, and Michel \cite{FKM} have already proved the square root cancellation in the character sum of type $\mathcal{C}'_{2, i}(...)$ when the tuple $(\gamma_1, \gamma_2, \gamma_3, \gamma_4)$ is normal, i.e., when the cardinality of the set $ S \, =\,\{1 \leq i \leq 4:\,\gamma_i \equiv \gamma\, (\bmod\, p_i)\}$ is odd for some $\gamma \in {PGL_2}(\mathbb{F}_p)$. There are four possible values of $|S|$. If $|S| = 1$ or $|S| = 3$ (i.e., all the $\gamma_i$'s are different or any three of them are equal), then the tuple $(\gamma_1, \gamma_2, \gamma_3, \gamma_4)$ will be normal. In these cases, using  of \cite[Corollary $1.6 (1)$]{FKM}, we obtain our desired result, i.e.,
\begin{align}\label{m11}
 \mathcal{C}'_{2,i}(...) \,\ll\, p^{5/2}_i.   
\end{align}
Now we consider the remaining two cases, i.e., when $|S| = 2$ or $|S| = 4$. In the case when $|S| = 2$, we have two sub-cases. In the first sub-case, without loss of generality (WLOG), we consider $\gamma_1 = \gamma_2$ and $\gamma_3 \neq \gamma_4$. In this sub-case, we will get the congruence condition $n''_1s_1 \equiv n''_2s_2(\mathrm{mod}\, p_i)$. Now, to find the bound of the character sum $\mathcal{C}'_{2, i}(...) $, we will take the average over $n''_1, n''_2, n''_3, r_2$ sums with the occurred congruence condition and use Weil's bound for Kloosterman sums. As we are saving $p_i$ from $n''_1$ or $n''_2$ length. So, in this sub-case, we are saving $p^{1/2}_i$ extra as compared to the bound in \eqref{m11}. Now for the second sub-case, again WLOG, we take $\gamma_1 = \gamma_2$ and $\gamma_3 = \gamma_4$. We will get two congruence conditions with modulus $p_i$. Again estimating the character sum on average as in the above sub-case, we will save $p_i$ extra as compared to the bound in \eqref{m11}. For the fourth and last case when $|S| = 4$, by proceeding similarly, we find that the bound of the character sum $\mathcal{C}'_{2, i}(...)$ given in \eqref{m11} will again dominate. This completes the proof of the lemma.
\end{proof}
\vspace{0.3cm}

\noindent Using the bound of character sum $\mathcal{C}'_{2, i}(...)$ given in Lemma \ref{t1} into equation \eqref{m2}, we get   
\begin{align*}
\mathcal{C}'_2(...) \, \ll \, (q'_3)^{5/2}.    
\end{align*}
And hence,
\begin{align}\label{m3}
\mathcal{C}_2(...) \,\ll_{\epsilon}\, (q'_3)^{5/2}\,(q''_3)^{3+\epsilon}.  
\end{align}
Finally, by using the estimates from \eqref{m5} and \eqref{m3} into equation \eqref{m4}, we get
\begin{align}\label{m7}
\mathcal{C}(...) \,\ll_{\epsilon}\,\,\frac{q^{5/2}}{\sqrt{n'_1 n'_2n'_3r_1}}\,\,(q_1q_2q''_3)^{1/2}\,q^{\epsilon}.      
\end{align}
We use this bound for character sum $\mathcal{C}(...)$ in the final estimation.
\vspace{0.3cm}

\subsection{Analysis of the Integral transform}\label{M27}
From equation \eqref{L11}, we have
\begin{align*}
 \mathcal{I}(...) = &\,\int_{0}^\infty W_1(x_1)\, e\left(\frac{N_1ux_1}{qQ} \pm  \frac{3(N_1(n'_1)^2 n''_1 x_1)^{1/3}}{q} \right)\notag\\
&\times \,\int_{0}^\infty W_2(x_2)\, e\left(-\frac{N_1ux_2}{qQ} \pm  \frac{3(N_1(n'_2)^2 n''_2 x_2)^{1/3}}{q} \right)\,\notag\\
&\times \,\int_{0}^\infty W_3(x_3)\,W\left(\frac{N_2x_3- N_1x_1}{M}\right) \,\, e\left(- \frac{N_2ux_3}{qQ} \pm  \frac{3(N_2 (n'_3)^2 n''_3 x_3)^{1/3}}{q} \right)\notag\\
&\times\,\int_{0}^\infty W_4(x_4)\, e\left( \frac{N_3ux_4}{qQ} \pm  \frac{3(N_3 r_1^2 r_2 x_4)^{1/3}}{q} \right)dx_1\,dx_2\,dx_3\,dx_4.   
\end{align*}
To find a bound for the integral transform $\mathcal{I}(...)$, we break our analysis into two cases, namely $q < \sqrt{N} |u|$ and $q > \sqrt{N}|u|$.\\

\textbf{Case 1:} When $q < \sqrt{N} |u|$. In this case, first we consider the $x_2$-integral say $\mathcal{J}_2$, which is given by
	\begin{align*}
\mathcal{J}_2 = \int_{0}^\infty W_2(x_2)\, e\left(-\frac{N_1ux_2}{qQ} \pm  \frac{3(N_1(n'_2)^2 n''_2 x_2)^{1/3}}{q} \right).
\end{align*}
Changing the variables $x_2 \rightarrow x^3_2$, we can write
\begin{align*}
\mathcal{J}_2 = \int_{0}^\infty W_2(x^3_2)\,3x^2_2\,\, e\left(P_2(x_2) \right),    
\end{align*}
where the phase function $P_2(x_2)$ is given by
\begin{align*}
P_2(x_2) =   -\frac{N_1ux^3_2}{qQ} \pm  \frac{3(N_1(n'_2)^2 n''_2)^{1/3}}{q} x_2. 
\end{align*}
Differentiating twice and for $q < \sqrt{N} |u|$, we note that
\begin{align*}
P''_2(x_2) \,\asymp\,\frac{N|u|}{qQ} \gg \,1.
\end{align*}
Using the second derivative bound (see Lemma \ref{db}), we get
\begin{align}\label{sd1}
\mathcal{J}_2 \,\ll\,\sqrt{\frac{qQ}{N |u|}}.    
\end{align}
Proceeding with the same steps for the $x_4$-integral say $\mathcal{J}_4$, we get
\begin{align}\label{sd2}
\mathcal{J}_4 \,\ll\,\sqrt{\frac{qQ}{N |u|}}.    
\end{align}
Now, consider the $x_1$ and $x_3$ integrals say $\mathcal{J}_{1,3}$, we have
\begin{align*}
\mathcal{J}_{1, 3} =& \int_{0}^{\infty}\,\int_{0}^{\infty} W_1(x_1)\,W_3(x_3)\,W\left(\frac{N_2 x_3- N_1x_1}{M}\right)\\
&\times\,e\left(\frac{N_1ux_1}{qQ} \pm  \frac{3(N_1(n'_1)^2 n''_1 x_1)^{1/3}}{q} \right)\,\,e\left(- \frac{N_2ux_3}{qQ} \pm  \frac{3(N_2 (n'_3)^2 n''_3 x_3)^{1/3}}{q} \right)\,dx_1\,dx_3 .   \end{align*}
Using the change of the variables  $x_1 = v_1$ and  ${(N_2x_3- N_1x_1)}/{M} = v_2$, the determinant of the corresponding Jacobian matrix becomes $M/N$. So, we get
\begin{align*}
\mathcal{J}_{1, 3} \,=\,& \frac{M}{N}\,\int_{0}^{\infty}\,\int_{0}^{\infty}\, W_1(v_1)\,W_3\left( \frac{N_1}{N_2}v_1 + \frac{M}{N_2}v_2 \right)\,W\left(v_2\right)\\
&\times\,\,
 e\left(\frac{N_1u v_1}{qQ} \pm  \frac{3(N_1(n'_1)^2 n''_1 v_1)^{1/3}}{q} \right)\\
 &\times\,e\left(- \frac{N_2u}{qQ}\,\left(\frac{N_1}{N_2}v_1 + \frac{M}{N_2}v_2 \right) \pm  \frac{3(N_2 (n'_3)^2 n''_3)^{1/3}}{q}\,\left(\frac{N_1}{N_2}v_1 + \frac{M}{N_2}v_2 \right)^{1/3}  \right)\,du\,dv.     
\end{align*}
Now, estimating the integral $\mathcal{I}(...)$ using \eqref{sd1}, \eqref{sd2} and $v_1, v_2$ integrals trivially, we get the following bound for the case $q < \sqrt{N} |u|$
\begin{align}\label{IM1}
\mathcal{I}(...)\, \ll\,\frac{M}{N}\,\frac{q\,Q}{N |u|}N^{\epsilon}. 
\end{align}

\textbf{Case 2:} When $q > \sqrt{N} |u|$. In this case, we first consider the $x_1$ and $x_3$-integral and use the same change of variable i.e. $v_1 = x_1$ and $v_2 = (N_2x_3 - N_1x_1)/M$. Now, estimating the $v_1$ and $v_2$ integrals trivially along with $x_2$ and $x_4$ integrals, we get
\begin{align}\label{IM2}
\mathcal{I}(...) \,\ll\, \frac{M}{N} N^{\epsilon}.    
\end{align}

\subsection{Final estimates}
From equation \eqref{m7}, we have
\begin{align*}
    \eta(N, M) = &\frac{N^{8/3}}{Q}\sum_{q\leq Q}\frac{1}{q^{5}}\,\int_{|u| \leq Q^{\epsilon}} \, \psi(q,u) \,\,  \sum_{\pm}\,\sum_{n'_{1}|q}\,(n'_{1})^{1/3}\, \sum_{(n'_1)^2n''_{1}\ll M'_0}\,  \frac{{A_{\pi}(n''_{1},n'_{1})}}{ (n''_{1})^{1/3}} \notag\\
 &\times \sum_{\pm}\,\sum_{n'_{2}|q}\,(n'_{2})^{1/3}\, \sum_{(n'_2)^2n''_{2} \ll N'_0}\,  \frac{\overline{A_{\pi}(n''_{2},n'_{2})}}{ (n''_{2})^{1/3}}    \,\sum_{\pm}\,\sum_{n'_{3}|q}\,(n'_{3})^{1/3}\, \sum_{(n'_3)^2n''_{3} \ll M'_0}\,  \frac{A_{\pi}(n''_{3},n'_{3})}{ (n''_{3})^{1/3}}\, \notag \\
 &\times \,\,\sum_{\pm} \sum_{r_{1}|q}\,r^{1/3}_{1} \, \sum_{r^2_1r_{2}\ll N'_0}  \frac{\overline{A_{\pi}(r_{2},r_{1})}}{ r^{1/3}_{2}}\,\, \mathcal{C}(...)\, \,\mathcal{I}(...)\,du + O(X^{-A}),
\end{align*}
where,
\begin{align*}
N'_0 = {\max}\left\{\frac{q^3}{N}, N^{1/2}|u|^3\right\}N^{\epsilon} = 
	\begin{cases}
		\frac{q^3}{N}N^{\epsilon} &\text{if}\,\,q > \sqrt{N}|u| \\
		\sqrt{N} |u|^3 N^{\epsilon}   &\text{if}\,\, q < \sqrt{N} |u|,
	\end{cases}
 \end{align*} 
 \begin{align*}M'_0 \,=\,{\max}\left\{\frac{q^3}{N}, \frac{q^3}{N}\left(\frac{N}{M}\right)^3,  N^{1/2}|u|^3\right\} N^{\epsilon} = \frac{q^3}{N}\left(\frac{N}{M}\right)^3 N^{\epsilon}.   
\end{align*} 
We will estimate $\eta(N, M)$ into the two cases, namely the generic size case $q > \sqrt{N}|u|$ and the small size case $q < \sqrt{N}|u|$ and denote the contribution of these cases to $\eta(N, M)$ by $\eta^{gen}(N, M)$ and $\eta^{small}(N, M)$, respectively. So, we have
\begin{align}\label{fbb}
  \eta(N, M) \,\leq\,\eta^{small}(N, M) + \eta^{gen}(N, M).   
\end{align}

By plugging the bounds for the character sum $\mathcal{C}(...)$ given in \eqref{m7} and the integral transform $\mathcal{I}(...)$ given in \eqref{IM1} into the above expression of $\eta(N, M)$, we get
 \begin{align*}
\eta^{small}(N, M) \ll\, &N^{\epsilon}\,\frac{N^{8/3}}{Q}\,\frac{QM}{N^{2}}\,\int_{|u| \leq Q^{\epsilon}} \, \frac{|\psi(q,u)|}{|u|} \,\sum_{q\leq Q}\frac{1}{q^{5}}\,\sum_{n'_{1}|q}\,\, \sum_{(n'_1)^2 n''_{1}\ll M'_0}\,  \frac{|{A_{\pi}(n''_{1},n'_{1})}|}{(n'_{1})^{-1/3} (n''_{1})^{1/3}} \notag\\
 &\times \,\sum_{n'_{2}|q}\,\, \sum_{(n'_2)^2n''_{2} \ll N'_0}\,  \frac{\overline{|A_{\pi}(n''_{2},n'_{2})}|}{(n'_{2})^{-1/3} (n''_{2})^{1/3}}    \,\,\sum_{n'_{3}|q}\,(n'_{3})^{1/3}\, \sum_{(n'_3)^2n''_{3} \ll M'_0}\,  \frac{|A_{\pi}(n''_{3},n'_{3})|}{ (n''_{3})^{1/3}}\, \notag \\
 &\times \,\, \sum_{r_{1}|q}\, \, \sum_{r^2_1r_{2}\ll N'_0}  \frac{\overline{|A_{\pi}(r_{2},r_{1})}|}{r^{-1/3}_{1} r^{1/3}_{2}}\,\, \frac{q^{5/2}}{\sqrt{n'_1 n'_2n'_3r_1}}\,\,(q_1q_2q''_3)^{1/2}\,\,q\,\,du.  
 \end{align*}
We note that
\begin{align}\label{n2}
\sum_{n'_{2}|q}\,&\sum_{(n'_2)^2 n''_{2} \ll N'_0}  \frac{|A_{\pi}(n''_2, n'_2)|}{(n'_{2})^{-1/3}(n''_{2})^{1/3}}\,\frac{1}{\sqrt{n'_2}} \, \leq \,\mathop{\sum_{n'_{2}}\,\sum_{n''_2}}_{(n'_2)^2n''_{2} \ll N'_0}\,\frac{|A_{\pi}(n''_2, n'_2)|}{(n'_{2})^{1/6}\, (n''_{2})^{1/3}}\notag\\
&\ll\,\left(\mathop{\sum_{n'_{2}}\, \sum_{n''_{2}}}_{(n'_2)^2n''_{2} \ll N'_0} |A_{\pi}(n''_2,n'_2)|^2\right)^{1/2}\,\left(\mathop{\sum_{n'_{2}\leq q}\, \sum_{n''_{2}}}_{(n'_2)^2n''_{2} \ll N'_0} \frac{1}{((n'_2)^{1/6}(n''_2)^{1/3})^{2}}\right)^{1/2} \ll\, (N'_0)^{2/3}.
\end{align}
The same bound is also true for $r_1, r_2$ sums. Also, following the same steps as above, we can get
\begin{align*}
\sum_{n'_{1}|q}\,\, \sum_{(n'_1)^2n''_{1}\ll M'_0}\,  \frac{|{A_{\pi}(n''_{1},n'_{1})}|}{(n'_{1})^{-1/3} (n''_{1})^{1/3}} \,\frac{1}{\sqrt{n'_1}}  \,\ll\, (M'_0)^{2/3}. 
\end{align*}
The same bound also holds for $n'_3, n''_3$-sums. Hence, $\eta^{small}(N, M)$ is dominated by
\begin{align*}
\,{N^{2/3+\epsilon}}\,M\,\int_{|u| \leq Q^{\epsilon}} \, \frac{|\psi(q,u)|}{|u|} \,\sum_{q\leq Q}\frac{(q_1q_2q''_3)^{1/2}}{q^{3/2}}\, \left(\sqrt{N} |u|^3\right)^{4/3}\,\left(\frac{q^3}{N}\left(\frac{N}{M}\right)^3\right)^{4/3}\,du.
\end{align*}
Using $q < \sqrt{N} |u|$ and  writing $q = q_0q_1q'_3$, where $q_0 = q_2q''_3$ is the square full part of $q$, we can write 
\begin{align*}
\eta^{small}(N, M)\,\ll&\,\,N^{\epsilon}\,\frac{N^{21/4}}{M^3}\,\int_{|u| \leq Q^{\epsilon}} \, {|\psi(q,u)}|\, |u|^{13/2}\,du\\
&\times\,\,\mathop{\sum_{q_{1}| [n'_1,n'_2,n'_3,r_1] }}\,q^{1/2}_1\,\mathop{\sum_{q'_3 \leq Q/q_1}}_{q'_3 \text{is square free}}\,\mathop{\sum_{q_0 \leq Q/q_1q'_3}}_{4q_0 \text{is square full}}\,q^{1/2}_0.
\end{align*}
Now, using the same steps as in Lemma \ref{ll2}, we can get
\begin{align*}
 \mathop{\sum_{q_0 \leq Q/q_1q'_3}}_{q_0\, \text{is square full}}\,q^{1/2}_0\,\ll\, \frac{Q}{q_1q'_3}\, Q^{\epsilon}.   
\end{align*}
Using this estimate, we get
\begin{align*}
\eta^{small}(N, M)\, \ll\,N^{\epsilon}\,\frac{N^{21/4}\,Q}{M^3}\,\int_{|u| \leq Q^{\epsilon}} \, {|\psi(q,u)}|\, |u|^{13/2}\,du.
\end{align*}
From equation \eqref{l2b}, we know that the $L^1$ norm of the function $\psi(q, x)$ is bounded by $Q^{\epsilon}$. Using $Q \asymp \sqrt{N}$, we finally get
\begin{align}\label{fb5}
\eta^{small}(N, M)\, \ll\,N^{\epsilon}\,\frac{N^{23/4}}{M^3}.    
\end{align}

Now we will estimate $\eta(N, M)$ in the generic size case $q > \sqrt{N} |u|^3$. Using the bound of the character sum from \eqref{m7} and the integral transform \eqref{IM2}, we get
\begin{align*}
  \eta^{gen}(N, M) \ll\, &N^{\epsilon}\,\frac{N^{8/3}}{Q}\,\frac{M}{N}\,\int_{|u| \leq Q^{\epsilon}} \,{|\psi(q,u)|} \,\sum_{q\leq Q}\frac{1}{q^{5}}\,\sum_{n'_{1}|q}\,\, \sum_{(n'_1)^2 n''_{1}\ll M'_0}\,  \frac{|{A_{\pi}(n''_{1},n'_{1})}|}{(n'_{1})^{-1/3} (n''_{1})^{1/3}} \notag\\
 &\times \,\sum_{n'_{2}|q}\,\, \sum_{(n'_2)^2n''_{2} \ll N'_0}\,  \frac{\overline{|A_{\pi}(n''_{2},n'_{2})}|}{(n'_{2})^{-1/3} (n''_{2})^{1/3}}    \,\,\sum_{n'_{3}|q}\,(n'_{3})^{1/3}\, \sum_{(n'_3)^2n''_{3} \ll M'_0}\,  \frac{|A_{\pi}(n''_{3},n'_{3})|}{ (n''_{3})^{1/3}}\, \notag \\
 &\times \,\, \sum_{r_{1}|q}\, \, \sum_{r^2_1r_{2}\ll N'_0}  \frac{\overline{|A_{\pi}(r_{2},r_{1})}|}{r^{-1/3}_{1} r^{1/3}_{2}}\,\, \frac{q^{5/2}}{\sqrt{n'_1 n'_2n'_3r_1}}\,\,(q_1q_2q''_3)^{1/2}\,\,du. 
 \end{align*}
Following \eqref{n2}, we can write
\begin{align*}
\eta^{gen}(N, M) &\ll\,N^{\epsilon}\frac{N^{8/3}}{Q}\,\frac{M}{N}\,\int_{|u| \leq Q^{\epsilon}} {|\psi(q,u)|} \sum_{q\leq Q}\frac{(q_1q_2q''_3)^{1/2}}{q^{5/2}}\left(\frac{q^3}{N}\right)^{4/3}\left(\frac{q^3}{N}\left(\frac{N}{M}\right)^3\right)^{4/3} du\\
&\ll\,N^{\epsilon}\,\frac{N^{3}}{M^3\,Q}\,\mathop{\sum_{q_{1}| [n'_1,n'_2,n'_3,r_1] }}\,\,q^6_1\,\mathop{\sum_{q'_3 \leq Q/q_1}}_{q'_3 \text{is square free}}\,(q'_3)^{11/2}\,\mathop{\sum_{q_0 \leq Q/q_1q'_3}}_{4q_0 \text{is square full}}\,q^6_0.
\end{align*}
We follow the same steps as in Lemma \ref{ll2} to get
\begin{align*}
\mathop{\sum_{q_0 \leq Q/q_1q'_3}}_{q_0\, \text{is square full}}\,q^6_0\,\ll\, \left(\frac{Q}{q_1q'_3}\right)^{13/2} .   
\end{align*}
Using these estimates with $Q \asymp \sqrt{N}$, we get
\begin{align}\label{fb51}
\eta^{gen}(N, M) &\ll\,N^{\epsilon}\, \frac{N^{23/4}}{M^3}.   
\end{align}
Using the bound of \eqref{fb5} and \eqref{fb51} into \eqref{fbb}, we get
\begin{align*}
\eta(N, M) \,\leq\, \eta^{small}(N, M) + \eta^{gen}(N, M) \,\ll\,N^{\epsilon}\, \frac{N^{23/4}}{M^3}.      
\end{align*}
Finally, plugging this bound of  $\eta(N, M)$ in the equation \eqref{L18}, we get
\begin{align*}
 S(N,M) \ll \, ||\alpha||_2\, ||\beta||_2\,M^{1/4}\, \left(\frac{N^{23/4}}{M^3} \right)^{1/4}. 
\end{align*}
The above bound of $S(N,M)$ is non-trivial when
\begin{align*}
 \,M^{1/4}\, \left(\frac{N^{23/4}}{M^3} \right)^{1/4} \ll\, (MN)^{1/2} \,\,\,\,\Longleftrightarrow \,\,\,\, N^{15/16} < M < N^{1-\epsilon}.   
\end{align*}
This proves our Theorem \ref{th1}.
\newpage

\section{Proof of Theorem \ref{th2}}

In Theorem \ref{th2}, we aim to prove the cancellations in the following sum
\begin{align*}
T(X, S) =  \sum_{X \leq n \leq 2X}\,\sum_{m \in S} \, \lambda(n)\, A_{\pi}(1,n+m),
\end{align*}
where $S \subset [X, 2X]$ is a sparse set such that  $|S| \asymp X^{1/k}$ for the real number $k >1$, $\lambda(n)$ is any complex sequence. By applying the Cauchy-Schwarz inequality to the sum over $n$, we get
\begin{align} \label{TXS}
\mathcal{T}_{X, \mathsf{S}} &\leqslant ||\lambda||_2\, \left({\sum_{n \in \mathbb{Z}} V_1\left({\frac{n}{X}}\right)\left|\sum_{m\in S} \mathcal{A}_f(n+m)\right|^2}\right)^{1/2} =   \, ||\lambda||_2\,\, T_1(X, S)^{1/2},
 \end{align}
where $V_1$ is positive bump function, supported on the interval $[-2,2]$ with $W_1(x) = 1$ for $x \in [-1,1]$, also satisfies $x^j W_1^{(j)}(x) \,\ll_j\,1$, for $j\geqslant 0$ and
\begin{align*}
 T_1(X, S)  =\,\sum_{n \in \mathbb{Z}}\,  \sum_{m_1 \in S}\,\sum_{m_2 \in S}\,\, A_{\pi}(1, n+m_1)\,\overline{A_{\pi}(1, n+m_2)} \,V_1\left(\frac{n}{X}\right).
\end{align*}
Changing the variables $n + m_2 \longrightarrow n$, we can write
\begin{align*}
T_1(X, S)  =  \sum_{n \ll X}\,\overline{A_{\pi}(1, n)}\,\,\sum_{m_1 \in S}\,\sum_{m_2 \in S}\,\, {A_{\pi}(1, n+m_1-m_2)}\,V_1\left(\frac{n-m_2}{X}\right).    
\end{align*} 
Applying the Cauchy-Schwarz inequality to the sum over $n$, we get that $ T_1(X, S)$ is dominated by
\begin{align*}
 \left(\sum_{n \ll X}\,\left|\overline{A_{\pi}(1, n)}\right|^2\right)^{1/2}\,   \left(\sum_{n \ll X}\,\left|\,\sum_{m_1 \in S}\,\sum_{m_2 \in S}\,{A_{\pi}(1, n+m_1-m_2)}\,V_1\left(\frac{n-m_2}{X}\right)\right|^2\right)^{1/2}.
\end{align*}
By using the Ramanujan bound on average from Lemma \ref{ramanubound}, we get
\begin{align*}
T_1(X, S) \ll&\, X^{1/2}\,\left(\sum_{n \ll X}\,\sum_{m_1 \in S}\,\sum_{m_2 \in S}\, {A_{\pi}(1, n+m_1-m_2)}\,V_1\left(\frac{n-m_2}{X}\right) \right.\\
& \hspace{3.5cm}\left.\, \times \sum_{m_3 \in S}\,\sum_{m_4 \in S}\, \overline{A_{\pi}(1, n+m_3-m_4)}\,V_2\left(\frac{n-m_4}{X}\right)\right)^{1/2}.      
\end{align*}
Since $V_1$ is a real-valued function, so $V_2 = V_1$. Here, we are using different notations for the copy of $V_1$ to be more precise about our steps. Changing the variable $n+m_3-m_4 \longrightarrow n$, we get
\begin{align*}
T_1(X, S) \ll&\, X^{1/2}\,\,  \left(\sum_{n \ll X}\,\overline{A(1, n)}\,\,\mathop{\sum_{m_i \in S}}_{1 \leq i \leq 2^2}\,{A_{\pi}(1, n+m_1-m_2-m_3+m_4)}\right.\\
& \hspace{5cm}\left.\times \,V_1\left(\frac{n-m_2-m_3+m_4}{X}\right)\,V_2\left(\frac{n-m_3}{X}\right) \right)^{1/2}.   
\end{align*}
We again apply the Cauchy-Schwarz inequality to the $n$-sum, use the Ramanujan bound on average for the GL$(3)$ Fourier coefficients $A_{\pi}(1,n)$, and change the variables, we get that $T_1(X, S)$ is dominated by the following expression 
\begin{align*}
X^{1/2 + 1/2^2}\,&\left(\sum_{n \ll X}\,\overline{A_{\pi}(1, n)}\,\,\mathop{\sum_{m_i \in S}}_{1 \leq i \leq 2^3}\,{A_{\pi}(1, n+m_1-m_2-m_3+m_4-m_5+m_6+m_7-m_8)} \right. \\
&\hspace{3.5cm} \left. \times\,V_1\left(\frac{n-m_2-m_3+m_4-m_5 +m_6+m_7-m_8}{X}\right) \,\right. \\
&\left. \hspace{3.5cm} \,\times\,V_2\left(\frac{n-m_3-m_5+m_6+m_7-m_8}{X}\right)\right.\\
& \hspace{3.5cm}\left.\,\times \,V_3\left(\frac{n-m_5}{X}\right)\,V_4\left(\frac{n-m_5+m_6-m_8}{X}\right)\right)^{1/2^2}.    
\end{align*}
\vspace{0.2cm}

\noindent In the argument of $A_{\pi}(1, n)$ in the above expression, we can see that half of $m_i$'s appear with negative signs and the other half with positive signs. For simplicity, we make the change of variables and note that $T_1(X, S)$ is dominated by the following expression
\vspace{0.3cm}
\begin{align*}
&X^{1/2 + 1/2^2}\\
&\times\, \left(\sum_{n \ll X}\,\overline{A_{\pi}(1, n)}\,\,\mathop{\sum_{m_i \in S}}_{1 \leq i \leq 2^3}\,{A_{\pi}(1, n+(m_1-m_2)+(m_3-m_4)+(m_5-m_6)+(m_7-m_8))}\right.\\
&  \left.\,\vspace{3cm}\,\times \,\prod_{1 \leq j \leq 2^2}\,V_j\left(\frac{n + \mathop{\sum}_{1 \leq s \leq 2^3}\,\epsilon_{js} m_s}{X}\right)\,\right)^{1/2^2},    
\end{align*}
where $\epsilon_{js} \in \left\{0, 1, -1\right\}$. Let $A = S-S: = \{s_1 - s_2 : s_1 \in S, s_2 \in S\} \subset [-X, X]$. We repeat the same treatment as above to the $n$-sum, $\ell$-many times such that $\ell$ is least with $ |2^{\ell} A| \gg X^{1-\wp}$ for some $\wp \geq 1/16$ and $|2^{\ell} A| \ll X^{2^{\ell+1}/k} < X^{1-\epsilon}$ for any epsilon $0<\epsilon < 1/16$. We finally get that $T_1(X, S)$ is dominated by
\begin{align*}
 &\,X^{1/2 + 1/2^2 + ... + 1/2^{\ell}}\\
 &\times\,\,\left(\sum_{{n \ll X}}\,\overline{A_{\pi}(1, n)}\,\,\mathop{\sum_{m_i \in S}}_{1 \leq i \leq 2^{\ell+1}}\,{A_{\pi}(1, n+(m_1-m_2)+ (m_3-m_4)+...+(m_{2^{\ell+1}-1}- m_{2^{\ell+1}}}) \right. \\
 & \left. \times \prod_{1 \leq j \leq 2^{\ell}}\,V_j\left(\frac{n + \mathop{\sum}_{1 \leq s \leq 2^{\ell +1}}\,\epsilon_{js} m_s}{X}\right)   \right)^{1/2^{\ell}}. 
 \end{align*}
 Here $V_j = V_1$ for all $1 \leq j \leq 2^{\ell}$. It is clear that we can explicitly write the expression at any stage $\ell$. Using the inverse Fourier transform for the function $V_j$ for all $1 \leq j \leq 2^{\ell}$, we get that $T_1(X, S)$ is dominated by
 \begin{align*}
&X^{1-1/2^\ell}
\Bigg(
\int_{|y_1|\ll X^\epsilon}
\int_{|y_2|\ll X^\epsilon}
\cdots
\int_{|y_{2^\ell}|\ll X^\epsilon}
\widehat{V_1}(y_1)\widehat{V_2}(y_2)
\cdots\widehat{V_{2^\ell}}(y_{2^\ell})
\\
&\qquad\times
\sum_{n\ll X}\overline{A_\pi(1,n)}
e\left(\frac{n(y_1+y_2+\cdots+y_{2^\ell})}{X}\right)
\\
&\qquad\times
\sum_{\substack{m_1,\ldots,m_{2^{\ell+1}}\in S}}
A_\pi\left(
1,n+(m_1-m_2)+(m_3-m_4)+\cdots+
(m_{2^{\ell+1}-1}-m_{2^{\ell+1}})
\right)
\\
&\qquad\times
\prod_{j=1}^{2^\ell}
e\left(
\frac{
\left(\sum_{s=1}^{2^{\ell+1}}\epsilon_{sj}m_s\right)y_j
}{X}
\right)
\Bigg)^{1/2^\ell}.
\end{align*}
where $\widehat{V_j}$ denotes the Fourier transform of the function $V_j$ for all $1 \leq j \leq 2^{\ell}$. Define
 \begin{align*}
 \alpha(n): = \overline{A_{\pi}(1,n)}\,e\left(\frac{n(y_1 + y_2 + ... + y_{2^{\ell})})}{X}\right),    
 \end{align*} 
 \begin{align*}
 \beta_{\ell}(m) :=  \mathop{\sum_{m_1 \in S}\,\sum_{m_2 \in S}\,...\,\sum_{m_{2^{\ell+1}} \in S}}_{m \,=\,(m_1-m_2) + (m_3 - m_4) \,+\, ...\, + \,(m_{2^{\ell+1}-1}- m_{2^{\ell+1}}) } \, \prod_{1 \leq j \leq 2^{\ell}}\,e\left(\frac{\left(\mathop{\sum}_{1 \leq s \leq 2^{\ell +1}}\,\epsilon_{sj} m_s\right)y_j}{X}\right). 
 \end{align*}
 Hence, we can write
\begin{align} \label{T1XS}
T_1(X, S) \ll \,X^{\epsilon}\,X^{1- 1/2^{\ell}}\, \left(\,\int_{|y_1| \ll X^{\epsilon}}\int_{|y_2| \ll X^{\epsilon}}...\int_{|y_{2^{\ell}}| \ll X^{\epsilon}}\,\widehat{V_1}(y_1)\,\widehat{V_2}(y_1)...\widehat{V_{2^{\ell}}}(y_{2^{\ell}})\,T_2(X, S)\right)^{1/2^{\ell}}, 
\end{align}
where $T_2(X, S) $ represents a bilinear sum given by
\begin{align} \label{T2XS}
T_2(X, S) = \,  \sum_{n \ll X}\,\sum_{m \leq |2^{\ell} A|}\, \alpha(n)\,\beta_{\ell}(m)\,A_{\pi}(1,n+m) . 
\end{align}
We note that the trivial bound for $T_2 (X, S)$ gives rise to 
\begin{align*}
     T_2 (X, S) \ll  X^\epsilon X \left|2^{\ell} A\right| \ll  X^\epsilon X \left|   A\right|^{2^{\ell}} \ll  X^\epsilon X \left|   S\right|^{2^{\ell+1}}. 
\end{align*}
Note that the Fourier transform $\widehat{V_j}$ is bounded as $V_j$ is a compactly supported smooth function for every $1 \leq j \leq 2^{\ell}$. Substituting the trivial bound for $  T_2 (X, S)$ in equation \eqref{T1XS} we obtain that
\begin{align*}
    T_1(X, S) \ll \,X^{1- 1/2^{\ell} + \epsilon}\, \left( X^{1+\epsilon} \left|   S\right|^{2^{\ell+1}}\right)^{1/2^{\ell}} \ll X^\epsilon X |S|^2. 
\end{align*}
Substituting the bound for $ T_1(X, S) $ in equation \eqref{TXS} we obtain that
\begin{align*}
    T(X, S) \ll ||\lambda||_2 \sqrt{X} |S| X^{\epsilon}. 
\end{align*}
We observe that a trivial bound for $ T_2(X, S)$ gives rise to a trivial bound for $T(X, S)$. We shall apply Theorem  \ref{th1} to get a non-trivial bound for  $T(X, S)$.
\vspace{0.3cm}

\noindent Since we have chosen $\ell$ (least) such that $|2^{\ell} A| \gg X^{1-\wp}$, i.e.
$X^{1-\wp} \ll |2^{\ell} A| \ll X^{2^{\ell+1}/k}$, which implies that $ 2^{\ell+1} > k(1-\wp)$. So, we can see that the above sum $T_2(X, S)$  is same as our main sum $S(N,M)$ in Theorem \ref{th1} with $N = X$ and $M = X^{2^{\ell+1}/k} < X^{1-\epsilon}$. By using Theorem \ref{th1} with $ ||\alpha||_2 \ll {X}^{1/2}$, we finally get 
\begin{align*}
 T(X, S) & \ll\, X^{\epsilon} ||\lambda||_2\, X^{1/2- 1/2^{\ell+1}} \,\left(||\alpha||_2\,||\beta_{\ell}||_2\, \frac{X^{15/16}}{X^{2^{\ell+1}/k}}\, \left(X\,X^{2^{\ell+1}/k}\right)^{1/2+\epsilon}\right)^{1/2^{\ell +1}} \\
 & \ll\, X^{\epsilon}  ||\lambda||_2\, \sqrt{X}\,  X^{\frac{15}{ 2^{\ell +5}}  }\,\left(\, \frac{||\beta_{\ell}||_2}{\left(X^{2^{\ell+1}/k}\right)^{1/2}} X^{\epsilon}\right)^{1/2^{\ell +1}}\\
 & \ll\,  ||\lambda||_2\,\, \sqrt{X}\,\,\left(||\beta_{\ell}||_2\right)^{1/2^{\ell+1}}\,  X^{\frac{15}{2^{\ell +5}}- \frac{1}{2k} + \epsilon}.
\end{align*}
This bound for $T(X, S)$ is non-trivial when
\begin{align*}
   \left(||\beta_{\ell}||_2\right)^{1/2^{\ell+1}}\,  X^{\frac{15}{2^{\ell +5}}- \frac{1}{2k} + \epsilon} < X^{\frac{1}{k}+ \epsilon}.
\end{align*}
This proves Theorem \ref{th2}.\\

{\bf Acknowledgement:} 
The authors would like to express their gratitude to Professor Ritabrata Munshi for discussing the problem and sharing his ideas. The authors also want to acknowledge the Indian Institute of Technology Kanpur for fostering an exceptional academic environment. The first author would also like to thank the Indian Institute of Information Technology Guwahati for providing excellent research facilities during the refinement of this manuscript. During the course of this research, S. K. Singh received partial support from the D.S.T. Inspire Faculty Fellowship under grant number DST/INSPIRE/04/2018/000945.

\end{document}